\documentclass[oneside,english]{amsart}
\usepackage{lmodern}
\usepackage[T1]{fontenc}
\usepackage[latin9]{inputenc}
\usepackage{geometry}
\usepackage{mathrsfs}
\usepackage{amsbsy}
\usepackage{amstext}
\usepackage{amsthm}
\usepackage{amssymb}
\usepackage{wasysym}

\makeatletter
\numberwithin{equation}{section}
\numberwithin{figure}{section}
\theoremstyle{plain}
\newtheorem{thm}{\protect\theoremname}
\theoremstyle{definition}
\newtheorem{defn}[thm]{\protect\definitionname}
\theoremstyle{remark}
\newtheorem{rem}[thm]{\protect\remarkname}
\theoremstyle{plain}
\newtheorem{prop}[thm]{\protect\propositionname}
\theoremstyle{plain}
\newtheorem{lem}[thm]{\protect\lemmaname}
\theoremstyle{plain}
\newtheorem{cor}[thm]{\protect\corollaryname}

\makeatother

\usepackage{babel}
\providecommand{\corollaryname}{Corollary}
\providecommand{\definitionname}{Definition}
\providecommand{\lemmaname}{Lemma}
\providecommand{\propositionname}{Proposition}
\providecommand{\remarkname}{Remark}
\providecommand{\theoremname}{Theorem}

\begin{document}
\global\long\def\e{e}%
\global\long\def\V{{\rm Vol}}%
\global\long\def\bs{\boldsymbol{\sigma}}%
\global\long\def\br{\boldsymbol{\rho}}%
\global\long\def\bp{\boldsymbol{\pi}}%
\global\long\def\btau{\boldsymbol{\tau}}%
\global\long\def\bx{\mathbf{x}}%
\global\long\def\by{\mathbf{y}}%
\global\long\def\bz{\mathbf{z}}%
\global\long\def\bv{\mathbf{v}}%
\global\long\def\bu{\mathbf{u}}%
\global\long\def\bi{\mathbf{i}}%
\global\long\def\bn{\mathbf{n}}%
\global\long\def\grad{\nabla_{sp}}%
\global\long\def\Hess{\nabla_{sp}^{2}}%
\global\long\def\lp{\Delta_{sp}}%
\global\long\def\gradE{\nabla_{\text{Euc}}}%
\global\long\def\HessE{\nabla_{\text{Euc}}^{2}}%
\global\long\def\HessEN{\hat{\nabla}_{\text{Euc}}^{2}}%
\global\long\def\ddq{\frac{d}{dR}}%
\global\long\def\qs{q_{\star}}%
\global\long\def\qss{q_{\star\star}}%
\global\long\def\lm{\lambda_{min}}%
\global\long\def\Es{E_{\star}}%
\global\long\def\EsN{E_{\star,N}}%
\global\long\def\EsNell{E_{\star,N_{\ell}}}%
\global\long\def\EH{E_{\Hess}}%
\global\long\def\Esh{\hat{E}_{\star}}%
\global\long\def\ds{d_{\star}}%
\global\long\def\Cs{\mathscr{C}_{\star}}%
\global\long\def\nh{\boldsymbol{\hat{\mathbf{n}}}}%
\global\long\def\BN{\mathbb{B}^{N}}%
\global\long\def\ii{\mathbf{i}}%
\global\long\def\SN{\mathbb{S}^{N-1}}%
\global\long\def\SM{\mathbb{S}^{M-1}}%
\global\long\def\SNq{\mathbb{S}^{N-1}(q)}%
\global\long\def\SNqd{\mathbb{S}^{N-1}(q_{d})}%
\global\long\def\SNqp{\mathbb{S}^{N-1}(q_{P})}%
\global\long\def\nd{\nu^{(\delta)}}%
\global\long\def\nz{\nu^{(0)}}%
\global\long\def\cls{c_{LS}}%
\global\long\def\qls{q_{LS}}%
\global\long\def\dls{\delta_{LS}}%
\global\long\def\E{\mathbb{E}}%
\global\long\def\P{\mathbb{P}}%
\global\long\def\R{\mathbb{R}}%
\global\long\def\spp{{\rm Supp}(\mu_{P})}%
\global\long\def\indic{\mathbf{1}}%
\global\long\def\lsc{\mu_{{\rm sc}}}%
\newcommand{\SNarg}[1]{\mathbb S^{N-1}(#1)}
\global\long\def\se{s(E)}%
\global\long\def\ses{s(\Es)}%
\global\long\def\so{s(0)}%
\global\long\def\sef{s(E_{f})}%
\global\long\def\seinf{s(E_{\infty})}%
\global\long\def\L{\mathcal{L}}%
\global\long\def\gflow#1#2{\varphi_{#2}(#1)}%
\global\long\def\S{\mathscr{S}}%
\global\long\def\Frep{F^{{\rm Rep}}}%
\global\long\def\s{\mathfrak{s}}%
\global\long\def\e{e}%

\title{TAP approach for multi-species spherical spin glasses I: general theory}
\author{Eliran Subag}
\begin{abstract}
We develop a generalized TAP approach for the multi-species version
of the spherical mixed $p$-spin models. In particular, we prove a
generalized TAP representation for the free energy at any overlap
vector which is multi-samplable in an appropriate sense. Moreover,
we show that if a multi-samplable overlap is maximal, then the TAP
correction is equal to an analogue of the well-known Onsager reaction
term. Finally, in a companion paper we use the results from the current
paper to compute the free energy at any temperature for all multi-species
pure $p$-spin models, assuming the free energy converges. 
\end{abstract}

\maketitle

\section{Introduction}

In the classical Sherrington-Kirkpatrick (SK) mean-field spin glass
model \cite{SK75}, the random interaction coefficients are identically
distributed for any two spins. Bipartite and, more generally, multi-species
versions of the SK model have been proposed in physics \cite{KorenblitShender85,KFS1,KFS2,Barra14,BarraGenoveseGuerra11,Barra}.
In those models the spins are divided into groups of different types,
and the strength of the interaction between any two spins depends
on their types. In this paper we will analyze the spherical mixed
$p$-spin version of those models defined as follows. 

First, consider a finite set of species $\S$, which will be fixed
throughout the paper. For each $N\geq1$,  suppose that
\[
\{1,\ldots,N\}=\bigcup_{s\in\S}I_{s},\quad\text{for some disjoint }I_{s}.
\]
The subsets $I_{s}$, of course, vary with $N$. Denoting $N_{s}=|I_{s}|$,
we will assume that the proportion of each species converges
\[
\lim_{N\to\infty}\frac{N_{s}}{N}=\lambda_{s}\in(0,1),\quad\text{for all }s\in\S.
\]

Let $S(d)=\{\bx\in\R^{d}:\,\|\bx\|=\sqrt{d}\}$ be the sphere of radius
$\sqrt{d}$ in dimension $d$. The \emph{configuration space} of the
spherical multi-species mixed $p$-spin model is
\[
S_{N}=\left\{ (\sigma_{1},\ldots,\sigma_{N})\in\R^{N}:\,\forall s\in\S,\,(\sigma_{i})_{i\in I_{s}}\in S(N_{s})\right\} .
\]
Denoting $\mathbb{Z}_{+}:=\{0,1,\ldots\}$ and $|p|:=\sum_{s\in\S}p(s)$
for $p\in\mathbb{Z}_{+}^{\S}$, let
\[
P=\Big\{ p\in\mathbb{Z}_{+}^{\S}:\,|p|\geq1\Big\}.
\]
Given some nonnegative numbers $(\Delta_{p})_{p\in P}$, define the
\emph{mixture polynomial} in the variables $x=(x(s))_{s\in\S}\in\R^{\S}$,
\[
\xi(x)=\sum_{p\in P}\Delta_{p}^{2}\prod_{s\in\S}x(s)^{p(s)}.
\]
We will assume that $\xi(1+\epsilon)<\infty$ for some $\epsilon>0$,
where for $c\in\R$ we write $\xi(c)$ for the evaluation of $\xi$
at the constant function $x\equiv c$. In analogy to the single-species
$p$-spin models, (i.e., with $|\S|=1$) we call the models such that
$\Delta_{p}^{2}>0$ for exactly one $p\in P$ multi-species \emph{pure}
$p$-spin models. 

We define the multi-species mixed \emph{$p$}-spin\emph{ Hamiltonian}
$H_{N}:S_{N}\to\R$ corresponding to the mixture $\xi$ by
\[
H_{N}(\bs)=\sqrt{N}\sum_{k=1}^{\infty}\sum_{i_{1},\dots,i_{k}=1}^{N}\Delta_{i_{1},\dots,i_{k}}J_{i_{1},\dots,i_{k}}\sigma_{i_{1}}\cdots\sigma_{i_{k}},
\]
where $J_{i_{1},\dots,i_{k}}$ are i.i.d. standard normal variables
and if $\#\{j\leq k:\,i_{j}\in I_{s}\}=p(s)$ for any $s\in\S$, then
$\Delta_{i_{1},\dots,i_{k}}$, which depends on $N$, is defined by
\begin{equation}
\Delta_{i_{1},\dots,i_{k}}^{2}=\Delta_{p}^{2}\frac{\prod_{s\in\S}p(s)!}{k!}\prod_{s\in\S}N_{s}^{-p(s)}.\label{eq:Delta}
\end{equation}
Here $\#A$ is the cardinality of a set $A$. By a straightforward
calculation, the covariance function of the random process $H_{N}(\bs)$
is given by
\begin{equation}
\frac{1}{N}\E H_{N}(\bs)H_{N}(\bs')=\xi(R(\bs,\bs')),\label{eq:cov}
\end{equation}
where we define the overlap vector
\[
R(\bs,\bs'):=\big(R_{s}(\bs,\bs')\big)_{s\in\S},\quad R_{s}(\bs,\bs'):=N_{s}^{-1}\sum_{i\in I_{s}}\sigma_{i}\sigma_{i}'.
\]

Identifying $S_{N}$ with the product space $\prod_{s\in\S}S(N_{s})$,
let $\mu$ be the product of the uniform measures on each of the spheres
$S(N_{s})$. The free energy $F_{N}$ and partition function $Z_{N}$
are defined by 
\begin{equation}
F_{N}:=\frac{1}{N}\log Z_{N}:=\frac{1}{N}\log\int_{S_{N}}e^{H_{N}(\bs)}d\mu(\bs).\label{eq:F}
\end{equation}
The Gibbs measure is the random probability measure on $S_{N}$ with
density
\[
\frac{dG_{N}}{d\mu}(\bs)=Z_{N}^{-1}e^{H_{N}(\bs)}.
\]

In the 70s, Thouless, Anderson and Palmer \cite{TAP} invented their
celebrated approach to analyze the Sherrington-Kirkpatrick model.
Their approach was further developed in physics, see e.g. \cite{TAP-SK1,TAP-SK3,TAP-pSPSG1,TAP-SK2,TAP-pS-Ising3,KurchanParisiVirasoro,Plefka},
with the general idea that for large $N$, $F_{N}(\beta)$ is approximated
by the free energies associated to the `physical' solutions of the
TAP equations. In particular, the single-species spherical pure $p$-spin
models were analyzed using the TAP approach non-rigorously by Kurchan,
Parisi and Virasoro in \cite{KurchanParisiVirasoro} and Crisanti
and Sommers in \cite{TAP-pSPSG1}. 

Recently, we introduced in \cite{FElandscape} a generalized TAP approach
for the single-species spherical models. The approach was extended
to mixed models with Ising spins by Chen, Panchenko and the author
\cite{TAPChenPanchenkoSubag,TAPIIChenPanchenkoSubag}. In the current
paper we will develop the generalized TAP approach for the multi-species
spherical models.

Analogously to the notion of multi-samplable overlaps introduced in
\cite{FElandscape} in the single-species setting, for the multi-species
models we use the following definition. Let $G_{N}^{\otimes n}$ denote
the $n$-fold product measure of $G_{N}$ with itself.
\begin{defn}
\label{def:mso}We will say that an overlap vector $q=(q(s))_{s\in\S}\in[0,1)^{\S}$
is multi-samplable if and only if for any $n\geq1$ and $\epsilon>0$,
\begin{equation}
\lim_{N\to\infty}\frac{1}{N}\log\E G_{N}^{\otimes n}\left\{ \forall i<j\leq n,\,s\in\S:\,\big|R_{s}(\bs^{i},\bs^{j})-q(s)\big|<\epsilon\right\} =0.\label{eq:multisamp-1}
\end{equation}
\end{defn}

Suppose that $q\in[0,1)^{\S}$ and define
\[
S_{N}(q)=\left\{ m\in\R^{N}:\,R(m,m)=q\right\} .
\]
Let $m\in S_{N}(q)$ be an arbitrary point. For any $\bs$ in 
\begin{equation}
\Big\{\bs\in S_{N}:\,\forall s\in\S,\,R_{s}(\bs,m)=R_{s}(m,m)\Big\},\label{eq:B0}
\end{equation}
define $\tilde{\bs}\in S_{N}$ by 
\begin{equation}
\tilde{\sigma}_{i}:=\sqrt{\frac{1}{1-q(s)}}(\sigma_{i}-m_{i}),\quad\text{if }i\in I_{s}.\label{eq:sigmatilde}
\end{equation}

Define the Hamiltonian 
\[
\tilde{H}_{N}(\tilde{\bs})=H_{N}(\bs)-H_{N}(m)
\]
on the image of (\ref{eq:B0}) by the mapping $\bs\mapsto\tilde{\bs}$,
which we may think of as a new multi-species model. Then, for any
two points $\bs^{1},\,\bs^{2}$ in (\ref{eq:B0}), by a straightforward
calculation,
\begin{equation}
\frac{1}{N}\E\tilde{H}_{N}(\tilde{\bs}^{1})\tilde{H}_{N}(\tilde{\bs}^{2})=\tilde{\xi}_{q}(R(\tilde{\bs}^{1},\tilde{\bs}^{2})),\label{eq:xiqcov}
\end{equation}
where we define the mixture 
\[
\tilde{\xi}_{q}(x)=\xi((1-q)x+q)-\xi(q)=\sum_{p\in P}\Delta_{q,p}^{2}\prod_{s\in\S}x(s)^{p(s)},
\]
where 
\[
\Delta_{q,p}^{2}:=\sum_{p'\geq p}\Delta_{p'}^{2}\prod_{s\in\S}\binom{p'(s)}{p(s)}(1-q(s))^{p(s)}q(s)^{p'(s)-p(s)},
\]
writing $p'\geq p$ if $p'(s)\geq p(s)$ for all $s\in\S$. Here,
operations between functions of $s$ are interpreted in the usual
way, namely, $\big((1-q)x+q\big)(s)=(1-q(s))x(s)+q(s)$. 

Additionally, define the mixture
\begin{equation}
\begin{aligned}\xi_{q}(x) & =\xi((1-q)x+q)-\xi(q)-\sum_{s\in\S}\big((1-q)\nabla\xi(q)x\big)(s)\\
 & =\sum_{p\in P:\,|p|\geq2}\Delta_{q,p}^{2}\prod_{s\in\S}x(s)^{p(s)},
\end{aligned}
\label{eq:xiq}
\end{equation}
where $\nabla\xi(q):=\Big(\frac{d}{dq(s)}\xi(q)\Big)_{s\in\S}$ and
the coefficients $\Delta_{q,p}^{2}$ are as above. 

With $H_{N}^{q}(\bs)$ being the Hamiltonian corresponding to $\xi_{q}$,
define the free energy
\[
F_{N}(q)=\frac{1}{N}\log\int_{S_{N}}e^{H_{N}^{q}(\bs)}d\mu(\bs).
\]
Lastly, define the ground state energy
\begin{equation}
\EsN(q):=\frac{1}{N}\max_{m\in S_{N}(q)}H_{N}(m).\label{eq:Esq}
\end{equation}
We will prove the following representation for the free energy.
\begin{thm}
[Generalized TAP representation]\label{thm:TAP} For any $q\in[0,1)^{\S}$,
$q$ is multi-samplable if and only if, as $N\to\infty$, 
\begin{equation}
\E F_{N}=\E\EsN(q)+\frac{1}{2}\sum_{s\in\S}\lambda_{s}\log(1-q(s))+\E F_{N}(q)+o(1).\label{eq:TAPrep}
\end{equation}
Moreover, for any $q\in[0,1)^{\S}$, 
\begin{equation}
\E F_{N}\geq\E\EsN(q)+\frac{1}{2}\sum_{s\in\S}\lambda_{s}\log(1-q(s))+\E F_{N}(q)+o(1).\label{eq:TAPrepineq}
\end{equation}
\end{thm}

Very recently, Bates and Sohn \cite{BatesSohn1,BatesSohn2} proved
a Parisi formula \cite{ParisiFormula,Parisi1980} for the limit of
the mean free energy $\E F_{N}$ for multi-species spherical mixed
$p$-spin models such that the mixture $\xi(x)$ is a convex function.
(For the proof of the formula in the single-species case, see \cite{Chen,Panch,Talag,Talag2}.)
The same convexity condition was required in an earlier work of Barra,
Contucci, Mingione and Tantari \cite{Barra}, where they proved for
the multi-species Sherrington-Kirkpatrick model that a Parisi type
formula is an upper bound for the limiting free energy, using an analogue
of Guerra's interpolation \cite{GuerraBound}. Panchenko proved the
matching lower bound in \cite{PanchenkoMulti} using an analogue of
the Aizenman-Sims-Starr scheme \cite{ASSscheme} (the proof of the
lower bound does not require $\xi(x)$ to be convex). He also proved
there that the overlaps $R_{s}(\bs,\bs')$ of independent samples
from the Gibbs measure are synchronized for different species. The
synchronization mechanism of \cite{PanchenkoMulti} was also used
in the spherical setting in \cite{BatesSohn1}. The combination of
\cite{Barra} and \cite{PanchenkoMulti} establishes the Parisi formula
for the multi-species Sherrington-Kirkpatrick model, assuming the
convexity of $\xi(x)$.

Denote the difference of the two sides of the TAP representation by
\[
\mathcal{D}_{N}=\E F_{N}-\Big(\E\EsN(q)+\frac{1}{2}\sum_{s\in\S}\log(1-q(s))+\E F_{N}(q)\Big).
\]
Given a mixture $\xi$, we will assume in the theorem below that 
\begin{equation}
\text{for any }q\in[0,1)^{\S},\text{ \ensuremath{\mathcal{D}_{N}} converges as \ensuremath{N\to\infty}.}\label{eq:conv}
\end{equation}
If $\xi(x)$ is convex on $[0,1]^{\S}$, one can check that $\xi_{q}(x)$
is also convex. In this case, the convergence follows from the results
of \cite{BatesSohn1}. For non-convex $\xi(x)$ the convergence was
not proved, but is certainly expected.

We will say that a multi-samplable $q\in[0,1)^{\S}$ is maximal if
any $q'\neq q$ such that $q'(s)\geq q(s)$ for all $s\in\S$ is not
multi-samplable. For maximal multi-samplable overlaps we prove that
the TAP correction is given by the following analogue of the Onsager
reaction term (see e.g. \cite{CrisantiSommersTAPpspin,KurchanParisiVirasoro}),
assuming the convergence (\ref{eq:conv}). In the proof of the theorem,
we will also show that $\xi_{q}$ is replica symmetric in the sense
that independent samples from the corresponding Gibbs measure typically
have overlap approximately zero. 
\begin{thm}
[Onsager correction] \label{thm:Onsager}Let $q\in[0,1)^{\S}$ be
a maximal multi-samplable overlap and assume the convergence (\ref{eq:conv}).
Then,

\begin{equation}
\E F_{N}(q)=\frac{1}{2}\xi_{q}(1)+o(1).\label{eq:onsager}
\end{equation}
\end{thm}

One of the main motivations of the current paper was to develop a
method that would allow one to compute the free energy for some class
of models which do not satisfy the convexity assumption on $\xi(x)$,
that was crucial to the analysis of \cite{Barra,BatesSohn1,PanchenkoMulti}.
In a companion paper \cite{TAPmulti2}  we compute the free energy
of the multi-species spherical \emph{pure} $p$-spin models at any
temperature, assuming convergence as in (\ref{eq:conv}), using the
TAP representation (\ref{eq:TAPrep}) and the explicit correction
from (\ref{eq:onsager}). One can easily check that for the pure models,
$\xi(x)$ is not convex anywhere in $[0,1]^{\S}$ (assuming that $p(s)\geq1$
for at least two species).

In the single-species case, the TAP representation from \cite{FElandscape}
was used in \cite{SubagTAPpspin} to compute the free energy of the
spherical pure $p$-spin models. For the mixed single-species spherical
models the free energy is given by the well-known Crisanti-Sommers
representation \cite{Crisanti1992} of Parisi formula \cite{Parisi,ParisiFormula},
first prove by Talagrand \cite{Talag} and extended by Chen \cite{Chen}.
There are several additional earlier results that are connected to
the TAP representation in the single-species case. Belius and Kistler
\cite{BeliusKistlerTAP} established a TAP representation for the
free energy for the spherical $2$-spin model with an external field.
In \cite{geometryGibbs} the free energy was computed for single-species
spherical pure $p$-spin models at low enough temperature, building
on results about the critical points \cite{A-BA-C,ABA2,2nd,pspinext},
and was shown to be given by the TAP representation. A similar result
was proved for mixed models close to pure by Ben Arous, Zeitouni and
the author \cite{geometryMixed}. For more recent results on the TAP
equations and representation for the single-species models with Ising
spins, see \cite{AuffingerJagannthSpinDist,AuffingerJagannathTAP,BolthausenTAP,BolthausenMorita,ChatterjeeTAP,ChenPanchenkoTAP,TalagrandBookI}.

Lastly, we mention earlier results recently proved for the multi-species
spherical models. Multi-species models with two species, i.e., $\S=\{s_{1},s_{2}\}$,
are called bipartite models. Baik and Lee \cite{BaikLeeBipartite}
computed the free energy and the limiting law of its fluctuations
for the pure $p$-spin spherical bipartite model with $p(s_{1})=p(s_{2})=1$
using tools from random matrix theory. Auffinger and Chen \cite{AuffingerChenBipartite}
proved that for mixed spherical bipartite models and in the presence
of an external field, if $\xi(1)$ and the strength of the field are
small enough, then the limiting free energy is given by the replica
symmetric solution of an analogue of the Crisanti-Sommers formula
\cite{Crisanti1992}. Certain bounds on the exponential growth rate
of the mean number of minimum points, known as `complexity', were
also derived in \cite{AuffingerChenBipartite}. McKenna \cite{McKennaComplexity}
computed the exact asymptotics for the mean complexity of minimum
and critical points. Kivimae \cite{Kivimae} proved that for large
energies, the second moment of the complexity matches its mean squared
at exponential scale, assuming $p(s_{1})$ and $p(s_{2})$ are large
enough. Combining \cite{McKennaComplexity} and \cite{Kivimae} yields
a variational formula for the ground state energy of the pure bipartite
models. An application of the second moment method to compute the
critical temperature of some multi-species models has been studied
in \cite{2ndmomentmulti}. As mentioned above, Bates and Sohn proved
a Parisi formula for the multi-species spherical mixed $p$-spin models
in \cite{BatesSohn1}, assuming the convexity of $\xi(x)$. In \cite{BatesSohn2}
they showed it admits a representation analogous to the Crisanti-Sommers
formula \cite{Crisanti1992}. 

In the setting of multi-species models with Ising spins, we also mentioned
above the results of Barra, Contucci, Mingione and Tantari \cite{Barra}
and Panchenko \cite{PanchenkoMulti} whose combination proves the
Parisi formula for the multi-species SK model, with the same convexity
assumption on $\xi(x)$. Bates, Sloman and Sohn \cite{BatesSlomanSohn}
study the replica symmetric phase of the latter formula. Mourrat derived
in \cite{MourratBipartite} an upper bound for the free energy of
the bipartite SK model with only inter-species interactions, whose
mixture is non-convex.

In the next section we state additional results related to the TAP
representation and outline the proof of the main results. In Section
\ref{sec:Concentration:-proof-of} we prove a concentration for the
TAP free energy, stated in Theorem \ref{thm:concentration} below.
Other results which we state in the next section will be proved in
Sections \ref{sec:pfTAPwoutcorr} and \ref{sec:comp_corr}. Finally,
in Section \ref{sec:pfMAin} we will prove Theorems \ref{thm:TAP}
and \ref{thm:Onsager}.

\section{\label{sec:2}Additional results and outline of the main proofs}

Denote the convex hull of $S_{N}$ by
\[
M_{N}=\Big\{(\sigma_{1},\ldots,\sigma_{N})\in\R^{N}:\,\sum_{i\in I_{s}}\sigma_{i}^{2}<N_{s},\,\forall s\in\S\Big\}.
\]
For any $m\in M_{N}$, real numbers $\delta,\,\rho\geq0$ and integer
$n\geq1$, define 
\begin{align*}
B(m,\delta) & =\left\{ \bs\in S_{N}:\,\forall s\in\S,\,\big|R_{s}(\bs,m)-R_{s}(m,m)\big|\leq\delta\right\} ,\\
B(m,n,\delta,\rho) & =\left\{ (\bs^{i})_{i\leq n}\in B(m,\delta)^{n}:\,\forall i\neq j,\,s\in\S,\,\big|R_{s}(\bs^{i},\bs^{j})-R_{s}(m,m)\big|\leq\rho\right\} .
\end{align*}
The set $B(m,\delta)$ is a product of spherical `bands', one for
each species. For $\delta=0$, it is equal to the set from (\ref{eq:B0}).
The elements of $B(m,n,\delta,\rho)$ are $n$-tuples of points from
$B(m,\delta)$, which are pairwise approximately orthogonal for each
species relative to $m$ for small $\delta$ and $\rho$, in the sense
that
\begin{equation}
\big|R_{s}(\bs^{i}-m,\bs^{j}-m)\big|\leq O(\delta+\rho).\label{eq:orthm}
\end{equation}

Using those sets, we associate two free energies to any $m$,
\begin{equation}
\begin{aligned}F_{N}(m,\delta) & =\frac{1}{N}\log\int_{B(m,\delta)}\exp\Big(H_{N}(\bs)-H_{N}(m)\Big)d\mu(\bs),\\
F_{N}(m,n,\delta,\rho) & =\frac{1}{Nn}\log\int_{B(m,n,\delta,\rho)}\exp\Big(\sum_{i=1}^{n}\big(H_{N}(\bs^{i})-H_{N}(m)\big)\Big)d\mu(\bs^{1})\cdots d\mu(\bs^{n}).
\end{aligned}
\label{eq:Fm}
\end{equation}
Noting that
\begin{equation}
\begin{aligned} & F_{N}(m,n,\delta,\rho)=F_{N}(m,\delta)\\
 & +\frac{1}{Nn}\log G_{N}^{\otimes n}\Big\{\big|R_{s}(\bs^{i},\bs^{j})-R_{s}(m,m)\big|\leq\rho,\,\forall s\in\S,\,i\neq j\,\Big|\,\bs^{i}\in B(m,\delta)\Big\},
\end{aligned}
\label{eq:penalty}
\end{equation}
the free energy $F_{N}(m,n,\delta,\rho)$ defined with $n$ approximately
orthogonal replicas can be thought of as the one with a single replica
$F_{N}(m,\delta)$, plus a penalty term related to how the conditional
Gibbs measure is `spread' on $B(m,\delta)$.
\begin{rem}
\label{rem:sub-additive}By definition $F_{N}(m,n,\delta,\rho)$ decreases
as we decrease $\delta$ and $\rho$. Since $n\mapsto nF_{N}(m,n,\delta,\rho)$
is sub-additive, we also have that $F_{N}(m,kn,\delta,\rho)\leq F_{N}(m,n,\delta,\rho)$
for any $k\geq1$.
\end{rem}

It is well-known that the free energy $F_{N}$ concentrates around
its mean, see e.g. \cite[Theorem 1.2]{PanchenkoBook}. For any non-random
$m$, $F_{N}(m,\delta)$ and $F_{N}(m,n,\delta,\rho)$ are free energies
defined on the spaces $B(m,\delta)$ and $B(m,n,\delta,\rho)$, and
therefore they concentrate around their mean. The concentration of
the ground state energy $\EsN(q)$ around its mean follows from the
classical inequality of Borell-TIS \cite{Borell,TIS}. We state those
results in the proposition below for later use.
\begin{prop}
[Concentration \cite{Borell,TIS,PanchenkoBook}] \label{prop:conc_nonuni}For
any $q\in[0,1)^{\S}$, $n\geq1$, $\delta,\,\rho>0$, $m\in M_{N}$
and $t>0$, with $X_{N}$ being either $F_{N}$, $F_{N}(m,\delta)$,
$F_{N}(m,n,\delta,\rho)$ or $\EsN(q)$,
\[
\P\big(|X_{N}-\E X_{N}|\geq t\big)\leq2\exp\left(-\frac{Nt^{2}}{16\xi(1)}\right).
\]
\end{prop}

While $F_{N}(m,\delta)$ concentrates around its mean for any fixed
non-random $m$, the maximal deviation $F_{N}(m,\delta)-\E F_{N}(m,\delta)$
over $S_{N}(q)$ is typically of order $O(1)$. In contrast to $F_{N}(m,\delta)$,
the other free energy we defined $F_{N}(m,n,\delta,\rho)$ satisfies
the following uniform concentration result, crucial to our analysis.
\begin{thm}
[Uniform concentration] \label{thm:concentration}For any $q\in[0,1)^{\S}$
and $t,\,c>0$, there exist $\rho_{0},\,\delta_{0}>0$ and $n_{0}\geq1$
such that for any $n\geq n_{0}$, $\rho\leq\rho_{0}$, and $\delta\leq\delta_{0}$
and any $N$,
\begin{equation}
\P\Big(\max_{m\in S_{N}(q)}\big|F_{N}(m,n,\delta,\rho)-\E F_{N}(m,n,\delta,\rho)\big|>t\Big)<4e^{-Nc}.\label{eq:uconc}
\end{equation}
\end{thm}

For two sequences of random variables that (may) depend on $n$, $\delta$
and $\rho$, we will write $A_{N}(n,\delta,\rho)\approx B_{N}(n,\delta,\rho)$
if for any $\epsilon>0$, 
\begin{equation}
\lim_{\substack{n\to\infty\\
\delta,\rho\to0
}
}\liminf_{N\to\infty}\P\left\{ \big|A_{N}(n,\delta,\rho)-B_{N}(n,\delta,\rho)\big|<\epsilon\right\} =1.\label{eq:apxnotation}
\end{equation}
We will write $A_{N}(n,\delta,\rho)\apprle B_{N}(n,\delta,\rho)$
if the same holds without the absolute value. The same notation will
also be used for non-random sequences, in which case it will mean
that the corresponding bound holds for large enough $N$ deterministically. 

By definition, for any $m\in M_{N}$,
\begin{equation}
\frac{1}{N}H_{N}(m)+F_{N}(m,n,\delta,\rho)\leq\frac{1}{N}H_{N}(m)+F_{N}(m,\delta)\leq F_{N}.\label{eq:ineq}
\end{equation}
We will see (in Lemma \ref{lem:char2-1}) that $q$ is multi-samplable
if and only if there exists some (random) $m_{\star}\in S_{N}(q)$
for which 
\begin{equation}
\frac{1}{N}H_{N}(m_{\star})+F_{N}(m_{\star},n,\delta,\rho)\approx\frac{1}{N}H_{N}(m_{\star})+F_{N}(m_{\star},\delta)\approx F_{N}.\label{eq:apx}
\end{equation}
Combining (\ref{eq:ineq}) and (\ref{eq:apx}) with Proposition \ref{prop:conc_nonuni}
and Theorem \ref{thm:concentration} and the fact that the expectation
in (\ref{eq:uconc}) is constant on $S_{N}(q)$, one easily sees that
\begin{equation}
\frac{1}{N}H_{N}(m_{\star})\approx\EsN(q).\label{eq:GS_energy}
\end{equation}

We will prove the following characterization for multi-samplable overlaps
in Section \ref{sec:pfTAPwoutcorr}. Note that the fact that multi-samplability
implies (\ref{eq:char2}) follows from (\ref{eq:apx}), (\ref{eq:GS_energy})
and Theorems \ref{prop:conc_nonuni} and \ref{thm:concentration}.
\begin{prop}
\label{prop:TAPwithoutcorrection}An overlap vector $q\in[0,1)^{\S}$
is multi-samplable if and only if for arbitrary $m\in S_{N}(q)$,
\begin{equation}
\E F_{N}\approx\E\EsN(q)+\E F_{N}(m,n,\delta,\rho).\label{eq:char2}
\end{equation}
Moreover, for any $q\in[0,1)^{\S}$ ,
\[
\E F_{N}\gtrsim\E\EsN(q)+\E F_{N}(m,n,\delta,\rho).
\]
\end{prop}

The proof of Theorem \ref{thm:TAP} will be completed by the following
proposition.
\begin{prop}
\label{prop:correction}For any $q\in[0,1)^{\S}$,
\begin{equation}
\E F_{N}(m,n,\delta,\rho)\approx\frac{1}{2}\sum_{s\in\S}\lambda_{s}\log(1-q(s))+\E F_{N}(q).\label{eq:correction}
\end{equation}
\end{prop}

Let $H_{N}^{q}(\bs)$ and $\tilde{H}_{N}^{q}(\bs)$ be the Hamiltonians
corresponding to the mixtures $\xi_{q}$ and $\tilde{\xi}_{q}$ which
we defined in the introduction. Define
\[
F_{N}^{q}(n,\rho)=\frac{1}{Nn}\log\int_{B(0,n,\rho)}\exp\Big(\sum_{i=1}^{n}H_{N}^{q}(\bs^{i})\Big)d\mu(\bs^{1})\cdots d\mu(\bs^{n})
\]
and define $\tilde{F}_{N}^{q}(n,\rho)$ similarly with $\tilde{H}_{N}^{q}(\bs)$,
where $B(0,n,\rho):=B(0,n,\delta,\rho)$ is independent of $\delta$. 

The proof of Proposition \ref{prop:correction} will consist of the
following three lemmas. 
\begin{lem}
\label{lem:transfertoSN}For any $q\in[0,1)^{\S}$,
\begin{equation}
\lim_{\substack{n\to\infty\\
\delta,\rho\to0
}
}\limsup_{N\to\infty}\Big|\E F_{N}(m,n,\delta,\rho)-\Big(\E\tilde{F}_{N}^{q}(n,\rho)+\frac{1}{2}\sum_{s\in\S}\lambda_{s}\log(1-q(s))\Big)\Big|=0.\label{eq:lem8}
\end{equation}
\end{lem}

\begin{lem}
\label{lem:remove1spin}For any $q\in[0,1)^{\S}$,
\[
\lim_{\substack{n\to\infty\\
\rho\to0
}
}\limsup_{N\to\infty}\big|\E F_{N}^{q}(n,\rho)-\E\tilde{F}_{N}^{q}(n,\rho)\big|=0.
\]
\end{lem}

\begin{lem}
\label{lem:manyto1replica}For any $q\in[0,1)^{\S}$, $\rho>0$ and
$n\geq1$,
\begin{equation}
\limsup_{N\to\infty}\big|\E F_{N}(q)-\E F_{N}^{q}(n,\rho)\big|=0.\label{eq:manyto1replica}
\end{equation}
\end{lem}

The free energies $F_{N}(m,n,\delta,\rho)$ and $\tilde{F}_{N}^{q}(n,\rho)$
are computed using replicas from $B(m,\delta)$ and $S_{N}$, respectively.
Lemma \ref{lem:transfertoSN} will follow by mapping $B(m,\delta)$
to $B(m,0)$ and rescaling the latter to $S_{N-|\S|}$ in order to
relate the two free energies, and using the Lipschitz property of
the Hamiltonian to bound their difference. The logarithmic (entropy)
term in (\ref{eq:lem8}) will arise from a change of volumes, as the
limit
\[
\lim_{\delta\to0}\lim_{N\to\infty}\frac{1}{N}\log\mu(B(m,\delta))=\frac{1}{2}\sum_{s\in\S}\lambda_{s}\log(1-q(s)).
\]

\begin{proof}
[Proof of Lemma \ref{lem:remove1spin}] Of course, both Hamiltonians
$H_{N}^{q}(\bs)$ and $\tilde{H}_{N}^{q}(\bs)$ can be defined on
the same probability space such that 
\begin{equation}
\tilde{H}_{N}^{q}(\bs)=H_{N}^{q}(\bs)+\sum_{s\in\S}\sqrt{\frac{N}{N_{s}}}\Delta_{q,p_{s}}\sum_{i\in I_{s}}J_{i}\sigma_{i},\label{eq:1spindecom}
\end{equation}
where $J_{i}$ are i.i.d. standard normal variables independent of
$H_{N}^{q}(\bs)$ and $p_{s}$ is defined by $p_{s}(s')=1$ if $s'=s$
and otherwise $p_{s}(s')=0$. The point of this representation is
that as we let $n\to\infty$, the effect of the extra term (linear
in $\bs$) which is added to $H_{N}^{q}(\bs)$ becomes negligible
when computing $\tilde{F}_{N}^{q}(n,\rho)$. Precisely, for $(\bs^{1},\ldots,\bs^{n})$,
define $\bs:=\sum_{k=1}^{n}\bs^{k}$ and $J:=(J_{1},\ldots,J_{N})$.
Then, by Cauchy\textendash Schwarz, uniformly on $B(0,n,\rho)$,
\[
\begin{aligned}\sum_{k=1}^{n}\sum_{s\in\S}\sqrt{\frac{N}{N_{s}}}\Delta_{q,p_{s}}\sum_{i\in I_{s}}J_{i}\sigma_{i}^{k} & \leq\sum_{s\in\S}\Delta_{q,p_{s}}\sqrt{NN_{s}R_{s}(J,J)R_{s}(\bs,\bs)}\\
 & \leq\sqrt{N}\|J\|_{2}\sum_{s\in\S}\Delta_{q,p_{s}}\sqrt{R_{s}(\bs,\bs)},
\end{aligned}
\]
and 
\[
R_{s}(\bs,\bs)\leq n+n(n-1)\rho.
\]
fHence, almost surely,
\begin{equation}
\begin{aligned}\big|F_{N}^{q}(n,\rho)-\tilde{F}_{N}^{q}(n,\rho)\big| & \leq\end{aligned}
\sqrt{\left(\frac{1}{n}+\rho\right)\frac{1}{N}}\|J\|_{2}\sum_{s\in\S}\Delta_{q,p_{s}}.\label{eq:1sp}
\end{equation}
The proof is completed by Jensen's inequality since $\E\|J\|_{2}\leq(\E\|J\|_{2}^{2})^{1/2}=\sqrt{N}$.
\end{proof}
Lemma \ref{lem:remove1spin} allows us to work with the mixture $\xi_{q}$
instead of $\tilde{\xi}_{q}$, for which $\Delta_{p}=0$ if $|p|=1$.
For such mixtures we will prove the following lemma. 
\begin{lem}
\label{lem:0mso}Suppose that $\xi$ is a mixture which does not contain
single-spin interactions, namely, such that $\Delta_{p}=0$ if $|p|=1$.
Then, $q\equiv0$ is multi-samplable.
\end{lem}

Lemma \ref{lem:manyto1replica} follows from Lemma \ref{lem:0mso}.
\begin{proof}
[Proof of Lemma \ref{lem:manyto1replica}] The lemma will follow if
we show that for a mixture $\xi$ as in Lemma \ref{lem:0mso}, for
$m=0$,
\begin{equation}
\limsup_{N\to\infty}\big|\E F_{N}-\E F_{N}(0,n,\delta,\rho)\big|=0,\label{eq:Fm0}
\end{equation}
since applied to $\xi_{q}$ this exactly gives (\ref{eq:manyto1replica}).
For such $\xi$, by Lemma \ref{lem:0mso} and Definition \ref{def:mso},
for any fixed $t>0$, with probability that does not decay exponentially
fast as $N\to\infty$,
\[
\frac{1}{N}\log G_{N}^{\otimes n}\left\{ \forall i<j\leq n,\,s\in\S:\,\big|R_{s}(\bs^{i},\bs^{j})\big|<\rho\right\} >-t.
\]
By (\ref{eq:penalty}), since $B(0,\delta)=S_{N}$, for any $t>0$,
with such probability,
\[
\big|F_{N}-F_{N}(0,n,\delta,\rho)\big|<t,
\]
and (\ref{eq:Fm0}) follows from the concentration as in Proposition
\ref{prop:conc_nonuni}. 
\end{proof}
In Section \ref{sec:comp_corr} we will prove Lemmas \ref{lem:transfertoSN}
and \ref{lem:0mso}. This will complete the proof of Proposition \ref{prop:correction}.

\section{\label{sec:Concentration:-proof-of}Concentration: proof of theorem
\ref{thm:concentration}}

Let $q\in[0,1)^{\S}$ and let $m$ be an arbitrary point in $S_{N}(q)$.
For $\delta,\,\rho>0$ and $n\geq1$, consider the random field $(\bs^{1},\ldots,\bs^{n})\mapsto\sum_{i=1}^{n}\big(H_{N}(\bs^{i})-H_{N}(m)\big)$
on $B(m,n,\delta,\rho)$. Its variance at any point satisfies
\begin{align*}
 & \frac{1}{N}\E\Big\{\Big(\sum_{i=1}^{n}\big(H_{N}(\bs^{i})-H_{N}(m)\big)\Big)^{2}\Big\}\\
 & =\sum_{i=1}^{n}\sum_{j=1}^{n}\Big(\xi(R(\bs^{i},\bs^{j}))-\xi(R(\bs^{i},m))-\xi(R(\bs^{j},m))+\xi(R(m,m))\Big)\\
 & \leq n\xi(1)+n^{2}\sqrt{|\S|}(2\delta+\rho)\|\nabla\xi(1)\|=:C(n,\rho,\delta),
\end{align*}
where the inequality follows since $R(m,m)=q$, $\|R(\bs^{i},m)-q\|\leq\sqrt{|\S|}\delta$
and $\|R(\bs^{i},\bs^{j})-q\|\leq\sqrt{|\S|}\rho$ for $i\neq j$. 

$F_{N}(m,n,\delta,\rho)$ is equal to $\frac{1}{n}$ times the free
energy of the random field above. Hence, from the well-known concentration
of the free energy \cite[Theorem 1.2]{PanchenkoBook}, for any fixed
$m\in S_{N}(q)$,
\[
\P\Big(\big|F_{N}(m,n,\delta,\rho)-\E F_{N}(m,n,\delta,\rho)\big|>t\Big)<2e^{-\frac{t^{2}n^{2}}{4C(n,\rho,\delta)}N}.
\]
If $q(s)=0$ for any $s$, this completes the proof. From now on,
assume otherwise.

By a union bound, for any $t,\,c,\,a>0$ and any subset $A_{N}\subset S_{N}(q)$
of cardinality $|A|\leq e^{aN}$, if $n$ is large enough and $\rho$
and $\delta$ are small enough,
\begin{equation}
\P\Big(\max_{m\in A_{N}}\big|F_{N}(m,n,\delta,\rho)-\E F_{N}(m,n,\delta,\rho)\big|>\frac{t}{2}\Big)<2e^{-2cN}.\label{eq:union}
\end{equation}
For any $m,\,m'\in S_{N}(q)$, $B(m,\delta)$ can be mapped to $B(m',\delta)$
by a mapping $\Theta:S_{N}\to S_{N}$ acting on the coordinates $(\sigma_{i})_{i\in I_{s}}$
as a rotation if $q(s)>0$ and as the identity if $q(s)=0$, such
that $\Theta(m)=m'$ and
\[
q(s)R_{s}(\bs-\Theta(\bs),\bs-\Theta(\bs))\leq R_{s}(m-m',m-m').
\]
In the Appendix we will prove a bound on the Lipschitz constant of
the Hamiltonian $H_{N}(\bs)$, see Lemma \ref{lem:Lipschitz}. From
this bound it follows that, for large enough $L>0$,
\begin{equation}
\begin{aligned} & \P\Big\{\forall m,m'\in S_{N}(q):\\
 & \,\,\,|F_{N}(m,n,\delta,\rho)-F_{N}(m',n,\delta,\rho)|\leq L\cdot D(m,m')\Big\}>1-e^{-2cN},
\end{aligned}
\label{eq:LF}
\end{equation}
where 
\[
D(m,m'):=\max_{s\in\S:q(s)>0}\sqrt{R_{s}(m-m',m-m')/q(s)}.
\]

Since $S_{N}(q)$ is a product of spheres, for any $x>0$ there exists
a subset $A_{N}\subset S_{N}(q)$ such that $|A_{N}|\leq e^{aN}$
for some large enough $a=a(x)$ and such that for any $m'\in S_{N}(q)$,
for some $m\in A_{N}$, $D(m,m')\leq x$. Assuming $A_{N}$ is such
a subset with $x=t/2L$, the proof is completed by (\ref{eq:union})
and (\ref{eq:LF}), since $\E F_{N}(m,n,\delta,\rho)$ is constant
on $S_{N}(q)$.\qed

\section{\label{sec:pfTAPwoutcorr}First characterization: Proof of Proposition
\ref{prop:TAPwithoutcorrection}}

In this section we prove Proposition \ref{prop:TAPwithoutcorrection}.
We will need the following lemma. 
\begin{lem}
\label{lem:char2}For any $q\in[0,1)^{\S}$, \textup{$q$ }is multi-samplable
if and only if for any $n\geq1$ and $\delta,\,\rho,\,t>0$, 
\begin{equation}
\lim_{N\to\infty}\frac{1}{N}\log\P\left\{ \exists m\in S_{N}(q):\,\Big|\frac{1}{N}H_{N}(m)+F_{N}(m,n,\delta,\rho)-F_{N}\Big|<t\right\} =0.\label{eq:goodm}
\end{equation}
\end{lem}

\begin{proof}
Let $q\in[0,1)^{\S}$ and assume that (\ref{eq:goodm}) holds for
any $n\geq1$ and $\delta,\,\rho,\,t>0$. We will prove that $q$
is multi-samplable. Note that 
\begin{equation}
\frac{1}{N}H_{N}(m)+F_{N}(m,n,\delta,\rho)-F_{N}=\frac{1}{nN}\log G_{N}^{\otimes n}\left\{ (\bs^{1},\ldots,\bs^{n})\in B(m,n,\delta,\rho)\right\} \leq0.\label{eq:FNGn}
\end{equation}
Hence, (\ref{eq:goodm}) is equivalent to
\begin{equation}
\lim_{N\to\infty}\frac{1}{N}\log\P\left\{ \sup_{m\in S_{N}(q)}G_{N}^{\otimes n}\left\{ (\bs^{1},\ldots,\bs^{n})\in B(m,n,\delta,\rho)\right\} >e^{-tnN}\right\} =0.\label{eq:Gn}
\end{equation}

By definition, if $(\bs^{1},\ldots,\bs^{n})\in B(m,n,\delta,\rho)$,
then for any $i<j\leq n$ and $s\in\S$,
\begin{align*}
|R_{s}(\bs^{i},\bs^{j})-q(s)| & \leq\rho.
\end{align*}
Hence,
\begin{align*}
 & \lim_{N\to\infty}\frac{1}{N}\log\E G_{N}^{\otimes n}\left\{ \forall i<j\leq n,\,s\in\S:\,\big|R_{s}(\bs^{i},\bs^{j})-q(s)\big|\leq\rho\right\} \\
 & \geq\lim_{N\to\infty}\frac{1}{N}\log\E\sup_{m\in S_{N}(q)}G_{N}^{\otimes n}\left\{ (\bs^{1},\ldots,\bs^{n})\in B(m,n,\delta,\rho)\right\} =0,
\end{align*}
where the equality follows from (\ref{eq:Gn}), and thus $q$ is multi-samplable.

Next, we will show that if $q$ is multi-samplable, then (\ref{eq:goodm})
holds. Let $q\in[0,1)^{\S}$, $n\geq1$ and $\epsilon>0$. Let $\bs^{1},\ldots,\bs^{2n}\in S_{N}$
be independent samples from $G_{N}$. Suppose that for any $i<j\leq2n$
and $s\in\S$,
\begin{equation}
\big|R_{s}(\bs^{i},\bs^{j})-q(s)\big|<\epsilon.\label{eq:Rq}
\end{equation}
Then $m':=\frac{1}{n}\sum_{i=n+1}^{2n}\bs^{i}$ satisfies
\[
\big|R_{s}(m',m')-q(s)\big|\leq\frac{1}{n}+\frac{n-1}{n}\epsilon,
\]
and for any $i<j\leq n$ and $s\in\S$,
\[
\big|R_{s}(\bs^{i},m')-q(s)\big|\leq\frac{1}{n}+\frac{n-1}{n}\epsilon
\]
and
\[
\big|R_{s}(\bs^{i}-m',\bs^{j}-m')-q(s)\big|<\epsilon+3\Big(\frac{1}{n}+\frac{n-1}{n}\epsilon\Big).
\]

If we define $m_{\star}=(m_{\star,i})_{i\leq N}$ by
\begin{equation}
m_{\star,i}:=\sqrt{\frac{q(s)}{R_{s}(m',m')}}m_{i}',\quad\text{if }i\in I_{s},\label{eq:mstar}
\end{equation}
(note that $R_{s}(m',m')>0$ a.s.) then $R_{s}(m_{\star},m_{\star})=q(s)$
and by Cauchy\textendash Schwarz,
\begin{align*}
|R_{s}(\bs^{i},m_{\star})-R_{s}(\bs^{i},m')| & \leq\sqrt{R_{s}(m_{\star}-m',m_{\star}-m')},\\
|R_{s}(m_{\star},m_{\star})-R_{s}(m',m')| & \leq2\sqrt{R_{s}(m_{\star}-m',m_{\star}-m')},
\end{align*}
and
\begin{align*}
R_{s}(m_{\star}-m',m_{\star}-m') & =\Big(\sqrt{R_{s}(m_{\star},m_{\star})}-\sqrt{R_{s}(m',m')}\Big)^{2}\\
 & \leq|R_{s}(m',m')-q(s)|,
\end{align*}
where the equality follows from (\ref{eq:mstar}).

By combining the above, given some $\delta$ and $\rho$, there exist
$n_{0}=n_{0}(\delta,\rho)$ and $\epsilon_{0}=\epsilon_{0}(\delta,\rho)$
such that for any $n\geq n_{0}$ and $\epsilon\leq\epsilon_{0}$ the
following holds. On the event that $\bs^{1},\ldots,\bs^{2n}\in S_{N}$
satisfy (\ref{eq:Rq}), for $m_{\star}$ as defined above using the
samples $\bs^{n+1},\ldots,\bs^{2n}$, (which are independent of $\bs^{1},\ldots,\bs^{n}$)
\[
m_{\star}\in S_{N}(q)\quad\text{and}\quad(\bs^{1},\ldots,\bs^{n})\in B(m_{\star},n,\delta,\rho).
\]
Therefore,
\begin{align*}
 & G_{N}^{\otimes2n}\left\{ \forall i<j\leq2n,\,s\in\S:\,\big|R_{s}(\bs^{i},\bs^{j})-q(s)\big|<\epsilon\right\} \\
 & \leq G_{N}^{\otimes2n}\left\{ (\bs^{1},\ldots,\bs^{n})\in B(m_{\star},n,\delta,\rho)\right\} \\
 & \leq\sup_{m\in S_{N}(q)}G_{N}^{\otimes n}\left\{ (\bs^{1},\ldots,\bs^{n})\in B(m,n,\delta,\rho)\right\} .
\end{align*}

Now, assume that $q$ is multi-samplable, let $\delta$ and $\rho$
be some positive numbers, and assume that $n\geq n_{0}(\delta,\rho)$
and $\epsilon\leq\epsilon_{0}(\delta,\rho)$. Then, for any $t>0$,
\[
\lim_{N\to\infty}\frac{1}{N}\log\P\left(G_{N}^{\otimes2n}\left\{ \forall i<j\leq2n,\,s\in\S:\,\big|R_{s}(\bs^{i},\bs^{j})-q(s)\big|<\epsilon\right\} >e^{-tN}\right)=0,
\]
and from the inequality above
\[
\lim_{N\to\infty}\frac{1}{N}\log\P\left(\exists m\in S_{N}(q):\,G_{N}^{\otimes n}\left\{ (\bs^{1},\ldots,\bs^{n})\in B(m,n,\delta,\rho)\right\} >e^{-tN}\right)=0.
\]

Combining this with (\ref{eq:FNGn}) implies (\ref{eq:goodm}), for
general $\delta$ and $\rho$ and $n\geq n_{0}(\delta,\rho)$. By
Remark \ref{rem:sub-additive}, (\ref{eq:goodm}) holds for general
$\delta$, $\rho$ and $n\geq1$.
\end{proof}
We continue with the proof of Proposition \ref{prop:TAPwithoutcorrection}.
Denote by $\mathcal{E}_{N}(t,n,\rho,\delta)$ the event 
\begin{align*}
\Big\{ & |F_{N}-\E F_{N}|<t,\,\Big|\EsN(q)-\E\EsN(q)\Big|<t,\\
 & \sup_{m\in S_{N}(q)}\Big|F_{N}(m,n,\delta,\rho)-\E F_{N}(m,n,\delta,\rho)\Big|<t\Big\}.
\end{align*}
Recall that by Proposition \ref{prop:conc_nonuni} and Theorem \ref{thm:concentration},
for any $t>0$ there exist $c$, $\delta_{0}$, $\rho_{0}$ and $n_{0}$,
such that if $\delta\leq\delta_{0}$, $\rho\leq\rho_{0}$ and $n\geq n_{0}$
and $N$ is large, then 
\begin{equation}
\P\left(\mathcal{E}_{N}(t,n,\rho,\delta)\right)\geq1-e^{-cN}.\label{eq:event_conc}
\end{equation}

Suppose that $m_{\star}\in S_{N}(q)$ is a point such that $\frac{1}{N}H_{N}(m_{\star})=\EsN(q)$.
From (\ref{eq:ineq}),
\[
F_{N}\geq\frac{1}{N}H_{N}(m_{\star})+F_{N}(m_{\star},n,\delta,\rho)=\EsN(q)+F_{N}(m_{\star},n,\delta,\rho).
\]
Combined with (\ref{eq:event_conc}), this implies that
\[
\E F_{N}\gtrsim\E\EsN(q)+\E F_{N}(m,n,\delta,\rho).
\]

Suppose that $q$ is multi-samplable. Let $t>0$ and suppose $c>0$,
$\delta_{0}$, $\rho_{0}$ and $n_{0}$ are the corresponding constants
from (\ref{eq:event_conc}). By Lemma \ref{lem:char2}, for any $\delta\leq\delta_{N}$,
$\rho\leq\rho_{0}$ and $n\geq n_{0}$, for large $N$, with positive
probability both $\mathcal{E}_{N}(t,n,\rho,\delta)$ occurs and there
exists $m\in S_{N}(q)$ such that 

\[
\Big|\frac{1}{N}H_{N}(m)+F_{N}(m,n,\delta,\rho)-F_{N}\Big|<t.
\]
Hence, 
\[
\E F_{N}\leq\E\EsN(q)+\E F_{N}(m,n,\delta,\rho)+4t
\]
and (\ref{eq:char2}) follows, since $t>0$ is arbitrary.

Conversely, if $q$ satisfies (\ref{eq:char2}), then on the event
$\mathcal{E}_{N}(t,n,\rho,\delta)$, 
\[
\Big|\frac{1}{N}H_{N}(m_{\star})+F_{N}(m_{\star},n,\delta,\rho)-F_{N}\Big|<4t,
\]
assuming that $n$ is large enough and $\rho$ and $\delta$ are small
enough, that $N$ is large, and that as before $m_{\star}\in S_{N}(q)$
is a point such that $\frac{1}{N}H_{N}(m_{\star})=\EsN(q)$. This
shows that given $t>0$, for $n$, $\rho$ and $\delta$ as before,
for large $N$,
\begin{equation}
\P\left\{ \exists m\in S_{N}(q):\,\Big|\frac{1}{N}H_{N}(m)+F_{N}(m,n,\delta,\rho)-F_{N}\Big|<4t\right\} \geq1-e^{-cN},\label{eq:goodm-1}
\end{equation}
where $c>0$ is a constant as in (\ref{eq:event_conc}). The fact
that the same holds for arbitrary $n$, $\delta$ and $\rho$ follows
from Remark \ref{rem:sub-additive} and (\ref{eq:ineq}). Thus, by
Lemma \ref{lem:char2}, $q$ is multi-samplable. This completes the
proof of Proposition \ref{prop:TAPwithoutcorrection}.\qed

We finish this section with the following improvement of Lemma \ref{lem:char2}.
\begin{lem}
\label{lem:char2-1}For any $q\in[0,1)^{\S}$, \textup{$q$ }is multi-samplable
if and only if for any $n\geq1$ and $\delta,\,\rho,\,t>0$, 
\begin{equation}
\limsup_{N\to\infty}\frac{1}{N}\log\P\left\{ \forall m\in S_{N}(q):\,\Big|\frac{1}{N}H_{N}(m)+F_{N}(m,n,\delta,\rho)-F_{N}\Big|>t\right\} <0.\label{eq:goodm-2}
\end{equation}
\end{lem}

\begin{proof}
If for some $q\in[0,1)^{\S}$ and any $n\geq1$ and $\delta,\,\rho,\,t>0$,
(\ref{eq:goodm-2}) holds, then (\ref{eq:goodm}) also holds and by
Lemma \ref{lem:char2}, $q$ is multi-samplable. 

Assume that $q\in[0,1)^{\S}$ is multi-samplable. In the proof above,
we have just seen that for such $q$, (\ref{eq:goodm-1}) holds for
arbitrary $n$, $\delta$, $\rho$ and $t$ and large $N$. From this,
(\ref{eq:goodm-2}) follows.
\end{proof}

\section{\label{sec:comp_corr}Computation of the correction: proof of Proposition
\ref{prop:correction}}

Proposition \ref{prop:correction} follows directly from Lemmas \ref{lem:transfertoSN}-\ref{lem:manyto1replica}
and the triangle inequality. In Section \ref{sec:2}, we already proved
Lemma \ref{lem:remove1spin} and proved that Lemma \ref{lem:manyto1replica}
follows from Lemma \ref{lem:0mso}. In this section we will prove
Lemmas \ref{lem:transfertoSN} and \ref{lem:0mso}, and by doing so
we will complete the proof of Proposition \ref{prop:correction}.

\subsection{An auxiliary lemma}

For the proof of Lemma \ref{lem:0mso}, we will need the following
lemma.
\begin{lem}
\label{lem:VarLB}Let $\theta>0$ and suppose that $\xi$ is a mixture
which does not contain single-spin interactions, namely, that $\Delta_{p}=0$
if $|p|=1$. If $\xi(1)<\theta$ and $N$ is such that $\max_{s}N/N_{s}<\theta$,
then for any $t>0$, 
\begin{equation}
\P\big(F_{N}\leq\E F_{N}-t\big)\leq e^{-C_{N}Nt^{2}},\label{eq:FNlb}
\end{equation}
where $C_{N}\to\infty$ is a sequence depending only on $\theta$.
\end{lem}

\begin{proof}
If $F_{N}^{k}$ is the free energy corresponding to the truncated
mixture $\xi^{k}(x)=\sum_{p:\,|p|\leq k}\Delta_{p}^{2}\prod_{s\in\S}x(s)^{p(s)}$,
then, $F_{N}^{k}\to F_{N}$ in probability and in $L^{1}$ as $k\to\infty$.
Therefore, it will be enough to the lemma assuming that $\Delta_{p}>0$
for infinitely many $p\in P$. 

Note that by H\"{o}lder's inequality, $F_{N}$ is a convex function
of the disorder variables $J_{i_{1},\dots,i_{k}}$. Hence, by the
main result of Paouris and Valettas \cite{paouris2018},
\[
\P\left(F_{N}\leq\E F_{N}-t\right)\leq\exp\left\{ -At^{2}/\text{Var}(F_{N})\right\} ,
\]
for some absolute constant $A>0$. Therefore, the lemma will follow
if we show that 
\[
\text{Var}(F_{N})\leq\frac{c_{N}}{N},
\]
for some sequence $c_{N}\to0$ depending only on $\theta$. In his
book \cite[Theorem 6.3]{ChattBook}, Chatterjee proved a bound of
this from (specifically, as in (\ref{eq:Vbound}) below) in the case
of the Sherrington-Kirkpatrick model. His method is general and in
the rest of the proof we adept it to our setting. For the usual spherical
models we used a similar proof in \cite{FElandscape}.

For any real (measurable) function $f\big((\bs^{i})_{i=1}^{n}\big)$
of $n$ replicas $\bs^{i}\in S_{N}$ denote by 
\[
\left\langle f\big((\bs^{i})_{i=1}^{n}\big)\right\rangle =\frac{\int_{S_{N}^{n}}f\big((\bs^{i})_{i=1}^{n}\big)e^{\sum_{i=1}^{n}H_{N}(\bs^{i})}d\bs^{1}\cdots d\bs^{n}}{\int_{S_{N}^{n}}e^{\sum_{i=1}^{n}H_{N}(\bs^{i})}d\bs^{1}\cdots d\bs^{n}}
\]
the average w.r.t. the Gibbs measure. 

Note that
\begin{equation}
\frac{\partial}{\partial J_{i_{1},\dots,i_{k}}}NF_{N}=\left\langle \sqrt{N}\Delta_{i_{1},\dots,i_{k}}\sigma_{i_{1}}\cdots\sigma_{i_{k}}\right\rangle .\label{eq:partialJ-1}
\end{equation}
For any function $f$ of $n$ replicas,
\begin{align*}
\frac{\partial}{\partial J_{i_{1},\dots,i_{k}}}\left\langle f\big((\bs^{a})_{a=1}^{n}\big)\right\rangle  & =\left\langle \sum_{a=1}^{n}\sqrt{N}\Delta_{i_{1},\dots,i_{k}}\sigma_{i_{1}}^{a}\cdots\sigma_{i_{k}}^{a}f\big((\bs^{a})_{a=1}^{n}\big)\right\rangle \\
 & -n\left\langle \sum_{a=1}^{n}\sqrt{N}\Delta_{i_{1},\dots,i_{k}}\sigma_{i_{1}}^{n+1}\cdots\sigma_{i_{k}}^{n+1}f\big((\bs^{a})_{a=1}^{n}\big)\right\rangle .
\end{align*}

Hence, by induction on $t\in\mathbb{Z}_{+}$, 
\begin{align*}
 & \frac{\partial}{\partial J_{i_{1}^{1},\dots,i_{k_{1}}^{1}}}\cdots\frac{\partial}{\partial J_{i_{1}^{t},\dots,i_{k_{t}}^{t}}}NF_{N}\\
 & =N^{\frac{t}{2}}\sum_{(\ell_{1},\ldots,\ell_{t})\in V(t)}c(\ell_{1},\ldots,\ell_{t})\left\langle \prod_{n=1}^{t}\Delta_{i_{1}^{n},\dots,i_{k_{n}}^{n}}\sigma_{i_{1}^{n}}^{\ell_{n}}\cdots\sigma_{i_{k_{n}}^{n}}^{\ell_{n}}\right\rangle ,
\end{align*}
where $V(t)=\{(\ell_{1},\ldots,\ell_{t})\in\mathbb{Z}^{t}:\,1\leq\ell_{n}\leq n\}$,
$c(\ell_{1},\ldots,\ell_{t})=\prod_{n=1}^{t}a(\ell_{n},n)$ and 
\[
a(\ell_{n},n)=\begin{cases}
1 & \text{if }\ell_{n}<n,\\
-(n-1) & \text{if }\ell_{n}=n.
\end{cases}
\]
In particular, $|c(\ell_{1},\ldots,\ell_{t})|\leq(t-1)!$.

Define 
\[
D_{t}:=\sum\Big(\E\frac{\partial}{\partial J_{i_{1}^{1},\dots,i_{k_{1}}^{1}}}\cdots\frac{\partial}{\partial J_{i_{1}^{t},\dots,i_{k_{t}}^{t}}}NF_{N}\Big)^{2}
\]
where the sum is over all possible choices of $k_{1},\ldots,k_{t}\geq1$
and $1\leq i_{1}^{n},\dots,i_{k_{n}}^{n}\leq N$.

Let $(\hat{\bs}^{i})_{i=1}^{t}$ be i.i.d. samples from the Gibbs
measure corresponding to an independent copy $\hat{H}_{N}(\hat{\bs})$
of the Hamiltonian $H_{N}(\bs)$. Then,
\begin{align*}
 & \Big(\sum_{(\ell_{1},\ldots,\ell_{t})\in V(t)}c(\ell_{1},\ldots,\ell_{t})\E\Big\langle\prod_{n=1}^{t}\sigma_{i_{1}^{n}}^{\ell_{n}}\cdots\sigma_{i_{k_{n}}^{n}}^{\ell_{n}}\Big\rangle\Big)^{2}\\
 & \leq(t!)^{3}\sum_{(\ell_{1},\ldots,\ell_{t})\in V(t)}\Big(\E\Big\langle\prod_{n=1}^{t}\sigma_{i_{1}^{n}}^{\ell_{n}}\cdots\sigma_{i_{k_{n}}^{n}}^{\ell_{n}}\Big\rangle\Big)^{2}\\
 & =(t!)^{3}\sum_{(\ell_{1},\ldots,\ell_{t})\in V(t)}\E\Big\langle\prod_{n=1}^{t}\sigma_{i_{1}^{n}}^{\ell_{n}}\hat{\sigma}_{i_{1}^{n}}^{\ell_{n}}\cdots\sigma_{i_{k_{n}}^{n}}^{\ell_{n}}\hat{\sigma}_{i_{k_{n}}^{n}}^{\ell_{n}}\Big\rangle.
\end{align*}

Combining the above, if we denote 
\[
I(k_{1},\ldots,k_{t}):=\bigg\{\left(i_{j}^{n}\right)_{\substack{1\leq n\leq t\\
1\leq j\leq k_{n}
}
}:\,i_{j}^{n}\in\{1,\ldots,N\},\,\forall j,n\bigg\},
\]
then
\begin{equation}
D_{t}\leq N^{t}(t!)^{3}\sum_{k_{1},\ldots,k_{t}}\sum_{(\ell_{1},\ldots,\ell_{t})\in V(t)}\sum_{I(k_{1},\ldots,k_{t})}\E\Big\langle\prod_{n=1}^{t}\Delta_{i_{1}^{n},\dots,i_{k_{n}}^{n}}^{2}\sigma_{i_{1}^{n}}^{\ell_{n}}\hat{\sigma}_{i_{1}^{n}}^{\ell_{n}}\cdots\sigma_{i_{k_{n}}^{n}}^{\ell_{n}}\hat{\sigma}_{i_{k_{n}}^{n}}^{\ell_{n}}\Big\rangle.\label{eq:Dtbd}
\end{equation}

Define 
\begin{align*}
I(p^{1},\ldots,p^{t}) & =\bigg\{\left(i_{j}^{n}\right)_{\substack{1\leq n\leq t\\
1\leq j\leq|p^{n}|
}
}:\,\#\{j:\,i_{j}^{n}\in I_{s}\}=p^{n}(s),\,\forall n,s\bigg\},\\
T(k_{1},\ldots,k_{t}) & =\Big\{(p^{1},\ldots,p^{t})\in P^{t}:\,|p^{n}|=k_{n}\Big\}.
\end{align*}
 Then
\begin{align*}
 & \sum_{I(k_{1},\ldots,k_{t})}\E\Big\langle\prod_{n=1}^{t}\Delta_{i_{1}^{n},\dots,i_{k_{n}}^{n}}^{2}\sigma_{i_{1}^{n}}^{\ell_{n}}\hat{\sigma}_{i_{1}^{n}}^{\ell_{n}}\cdots\sigma_{i_{k_{n}}^{n}}^{\ell_{n}}\hat{\sigma}_{i_{k_{n}}^{n}}^{\ell_{n}}\Big\rangle\\
 & =\sum_{T(k_{1},\ldots,k_{t})}\sum_{I(p^{1},\ldots,p^{t})}\E\Big\langle\prod_{n=1}^{t}\Delta_{i_{1}^{n},\dots,i_{|p^{n}|}^{n}}^{2}\sigma_{i_{1}^{n}}^{\ell_{n}}\hat{\sigma}_{i_{1}^{n}}^{\ell_{n}}\cdots\sigma_{i_{|p^{n}|}^{n}}^{\ell_{n}}\hat{\sigma}_{i_{|p^{n}|}^{n}}^{\ell_{n}}\Big\rangle.
\end{align*}

Since 
\[
\frac{\prod_{s\in\S}p(s)!}{|p|!}\prod_{s\in\S}N_{s}^{-p(s)}\prod_{(i_{j})_{j\leq|p|}\in I(p)}\sigma_{i_{1}}\hat{\sigma}_{i_{1}}\cdots\sigma_{i_{|p|}}\hat{\sigma}_{i_{|p|}}=\prod_{s\in\S}R_{s}(\bs,\hat{\bs})^{p(s)},
\]
we have that
\begin{align*}
 & \sum_{I(p^{1},\ldots,p^{t})}\E\Big\langle\prod_{n=1}^{t}\Delta_{i_{1}^{n},\dots,i_{|p^{n}|}^{n}}^{2}\sigma_{i_{1}^{n}}^{\ell_{n}}\hat{\sigma}_{i_{1}^{n}}^{\ell_{n}}\cdots\sigma_{i_{|p^{n}|}^{n}}^{\ell_{n}}\hat{\sigma}_{i_{|p^{n}|}^{n}}^{\ell_{n}}\Big\rangle\\
 & =\prod_{n=1}^{t}\Big(\Delta_{p^{n}}^{2}\frac{\prod_{s\in\S}p^{n}(s)!}{|p^{n}|!}\prod_{s\in\S}N_{s}^{-p^{n}(s)}\Big)\cdot\E\Big\langle\prod_{n=1}^{t}\sum_{(i_{j}^{n})_{j\leq|p^{n}|}\in I(p^{n})}\sigma_{i_{1}^{n}}^{\ell_{n}}\hat{\sigma}_{i_{1}^{n}}^{\ell_{n}}\cdots\sigma_{i_{|p^{n}|}^{n}}^{\ell_{n}}\hat{\sigma}_{i_{|p^{n}|}^{n}}^{\ell_{n}}\Big\rangle\\
 & =\Big(\prod_{n=1}^{t}\Delta_{p^{n}}^{2}\Big)\cdot\E\Big\langle\prod_{n=1}^{t}\prod_{s\in\S}R_{s}(\bs^{\ell_{n}},\hat{\bs}^{\ell_{n}})^{p^{n}(s)}\Big\rangle.
\end{align*}
And therefore,
\begin{align*}
 & \sum_{I(k_{1},\ldots,k_{t})}\E\Big\langle\prod_{n=1}^{t}\Delta_{i_{1}^{n},\dots,i_{k_{n}}^{n}}^{2}\sigma_{i_{1}^{n}}^{\ell_{n}}\hat{\sigma}_{i_{1}^{n}}^{\ell_{n}}\cdots\sigma_{i_{k_{n}}^{n}}^{\ell_{n}}\hat{\sigma}_{i_{k_{n}}^{n}}^{\ell_{n}}\Big\rangle\\
 & =\E\Big\langle\prod_{n=1}^{t}\sum_{p:|p|=k_{n}}\Delta_{p}^{2}\prod_{s\in\S}R_{s}(\bs^{\ell_{n}},\hat{\bs}^{\ell_{n}})^{p(s)}\Big\rangle.
\end{align*}

By combining this with (\ref{eq:Dtbd}), we obtain that

\begin{align*}
D_{t} & \leq N^{t}(t!)^{3}\sum_{k_{1},\ldots,k_{t}}\sum_{(\ell_{1},\ldots,\ell_{t})\in V(t)}\E\Big\langle\prod_{n=1}^{t}\sum_{p:|p|=k_{n}}\Delta_{p}^{2}\prod_{s\in\S}R_{s}(\bs^{\ell_{n}},\hat{\bs}^{\ell_{n}})^{p(s)}\Big\rangle.
\end{align*}

By H\"{o}lder's inequality, for $p^{1},\ldots,p^{t}\in P$ and $\gamma:=\sum_{n\leq t}\sum_{s\in\S}p^{n}(s)$,
\[
\E\Big\langle\prod_{n=1}^{t}\prod_{s\in\S}R_{s}(\bs^{\ell_{n}},\hat{\bs}^{\ell_{n}})^{p^{n}(s)}\Big\rangle\leq\max_{s\in\S}\E\Big\langle|R_{s}(\bs,\hat{\bs})|^{\gamma}\Big\rangle.
\]
Therefore, since that $|V(t)|=t!$,
\begin{align*}
D_{t} & \leq N^{t}(t!)^{4}\sum_{k_{1},\ldots,k_{t}}\Big(\max_{s\in\S}\E\Big\langle|R_{s}(\bs,\hat{\bs})|^{\sum_{n=1}^{t}k_{n}}\Big\rangle\prod_{n=1}^{t}\sum_{p:|p|=k_{n}}\Delta_{p}^{2}\Big)\\
 & \leq N^{t}(t!)^{4}\max_{s\in\S}\E\Big\langle|R_{s}(\bs,\hat{\bs})|^{2t}\Big\rangle\sum_{k_{1},\ldots,k_{t}}\prod_{n=1}^{t}\sum_{p:|p|=k_{n}}\Delta_{p}^{2}\\
 & =N^{t}(t!)^{4}\xi(1)^{t}\max_{s\in\S}\E\Big\langle|R_{s}(\bs,\hat{\bs})|^{2t}\Big\rangle,
\end{align*}
where in the second inequality we used the fact that $\Delta_{p}=0$
if $|p|=1$. 

Fix some $a\in(0,1)$. For large enough $N_{s}$, uniformly in $\gamma\leq N_{s}^{a}$,
for arbitrary $i\in I_{s}$,
\begin{align*}
\E\Big\langle|R_{s}(\bs,\hat{\bs})|^{\gamma}\Big\rangle & =N_{s}^{-\frac{\gamma}{2}}\E\Big\langle|\sigma_{i}|^{\gamma}\Big\rangle=N_{s}^{-\frac{\gamma}{2}}\frac{1}{\sqrt{N_{s}}}\frac{\omega_{N_{s}-1}}{\omega_{N_{s}}}\int_{-\sqrt{N_{s}}}^{\sqrt{N_{s}}}|x|^{\gamma}(1-x^{2}/N_{s})^{\frac{N_{s}-3}{2}}dx\\
 & \leq2N_{s}^{-\frac{\gamma}{2}}\int_{-\sqrt{N_{s}}}^{\sqrt{N_{s}}}|x|^{\gamma}e^{-\frac{1}{2}x^{2}}dx\leq2N_{s}^{-\frac{\gamma}{2}}\E|X|^{\gamma}\leq2N_{s}^{-\frac{\gamma}{2}}(\gamma-1)!!,
\end{align*}
where $\omega_{M}=2\pi^{M/2}/\Gamma(M/2)$ is the area of the unit
sphere in $\R^{M}$ and $X$ is a standard normal variable. Hence,
for appropriate constant $a>0$, for large $N$ and any $t\leq N^{a}$,
\begin{align*}
D_{t} & \leq4t^{6t}\theta^{2t}.
\end{align*}

From (\ref{eq:partialJ-1}), 
\begin{align*}
N^{2}\E|\nabla F_{N}|^{2} & :=\E\sum_{k}\sum_{1\leq i_{1},\dots,i_{k}\leq N}\Big(\frac{\partial}{\partial J_{i_{1},\dots,i_{k}}}NF_{N}\Big)^{2}\\
 & \leq N\E\sum_{k}\sum_{1\leq i_{1},\dots,i_{k}\leq N}\left\langle \Delta_{i_{1},\dots,i_{k}}^{2}\sigma_{i_{1}}^{2}\cdots\sigma_{i_{k}}^{2}\right\rangle =N\xi(1)\leq N\theta.
\end{align*}

By the bound \cite[Eq. (6.3)]{ChattBook} in Chatterjee's book, for
any $d\geq1$, 
\[
\text{Var}(NF_{N})\leq\sum_{t=1}^{d-1}\frac{D_{t}}{t!}+\frac{1}{d}N^{2}\E|\nabla F_{N}|^{2}.
\]
For $d=\left\lceil \epsilon\log N/\log\log N\right\rceil $, with
small enough $\epsilon>0$ depending on $\theta$, we obtain that
for some large $C>0$,
\begin{equation}
\text{Var}(NF_{N})\leq CN\frac{\log\log N}{\log N}.\label{eq:Vbound}
\end{equation}
\end{proof}

\subsection{Proof of Lemma \ref{lem:0mso}}

Let $n\geq1$ be an arbitrary integer. Suppose $\mathbf{X}=(\bx^{1},\ldots,\bx^{n})$
is a set of vectors in $S_{N}$ such that for any $i\neq j$ and $s\in\S$,
\begin{equation}
R_{s}(\bx^{i},\bx^{j})=0.\label{eq:R0}
\end{equation}
For $\delta>0$, define the set
\[
B(\delta,\mathbf{X})=\{\bs\in S_{N}:\,|R_{s}(\bs,\bx^{i})|\leq\delta,\,\forall s\in\S,\,i\leq n\}
\]
and the free energy
\[
F_{N}(\delta,\mathbf{X}):=\int_{B(\delta,\mathbf{X})}e^{H_{N}(\bs)}d\mu(\bs).
\]

We may identify $B(0,\mathbf{X})$ with a product of $|\S|$ spheres,
one for each $s\in\S$. Define $F_{N}(0,\mathbf{X})$ similarly to
the above with $\mu$ replaced by the product of the uniform measures
on those spheres. 

Obviously, there is a measure-preserving bijection that maps $B(0,\mathbf{X})$
to 
\begin{equation}
\left\{ (\sigma_{1},\ldots,\sigma_{N-n|\S|})\in\R^{N-n|\S|}:\,\forall s\in\S,\,(\sigma_{i})_{i\in I_{s}'}\in S(N_{s}-n)\right\} ,\label{eq:im}
\end{equation}
where $|I_{s}'|=|I_{s}|-n$ and $\cup_{s\in\S}I_{s}'=\{1,\ldots,N-n|\S|\}$,
which also preserves the overlap between any two points (defined on
the image space (\ref{eq:im}) using the subsets $I_{s}'$). Hence,
from (\ref{eq:cov}), $F_{N}(0,\mathbf{X})$ is equal to the free
energy of a multi-species Hamiltonian on (\ref{eq:im}) with mixture
$\xi$ multiplied by $\sqrt{\frac{N}{N-n|\S|}}$. (With this factor
accounting for the fact that the variance at any point is $N$. It
is easy to see that we may, and therefore will, neglect this factor
in our proof.) It follows by Lemma \ref{lem:VarLB} that for large
$N$ and any fixed $\mathbf{X}=(\bx^{1},\ldots,\bx^{n})$,
\[
\P\left\{ F_{N}(0,\mathbf{X})\leq\E F_{N}(0,\mathbf{X})-t\right\} \leq e^{-C_{N}Nt^{2}},
\]
where $C_{N}\to\infty$ is a sequence as in the same lemma. Of course,
\[
\lim_{N\to\infty}\big|\E F_{N}(0,\mathbf{X})-\E F_{N}\big|=0,
\]
and thus
\[
\P\left\{ F_{N}(0,\mathbf{X})\leq\E F_{N}-t\right\} \leq e^{-C_{N}Nt^{2}},
\]
where we may need to decrease $C_{N}$.

By a similar argument to the net argument used in the proof of Theorem
\ref{thm:concentration}, using the Lipschitz property of Lemma \ref{lem:Lipschitz}
and a union bound, one can see that for any fixed $t>0$ and $C>0$,
for large enough $N$,
\[
\P\left\{ \inf_{\mathbf{X}}(F_{N}(0,\mathbf{X})-\E F_{N})\leq-t\right\} \leq e^{-CN},
\]
where the infimum is over all $\mathbf{X}$ satisfying (\ref{eq:R0}).
Using the Lipschitz property of Lemma \ref{lem:Lipschitz} again,
it is easy to see that for any $t,\,C>0$, if $\delta>0$ is small
enough then for large $N$,
\[
\P\left\{ \inf_{\mathbf{X}}(F_{N}(\delta,\mathbf{X})-\E F_{N})\leq-t\right\} \leq e^{-CN}.
\]
From the concentration of $F_{N}$, decreasing the constant $C>0$
if needed, for large $N$,
\begin{equation}
\P\left\{ \inf_{\mathbf{X}}(F_{N}(\delta,\mathbf{X})-F_{N})\leq-t\right\} \leq e^{-CN}.\label{eq:Fframe}
\end{equation}

On the complement of the event in (\ref{eq:Fframe}), for any $k\leq n-1$,
\[
G_{N}^{\otimes k+1}\left\{ \max_{i\leq k,s\in\S}|R_{s}(\bs^{k+1},\bs^{i})|\leq\delta\,\Big|\,\bs^{1},\ldots\bs^{k}\right\} >e^{-tN}.
\]
Hence, with probability at least $1-e^{-CN}$ for large $N$,
\[
G_{N}^{\otimes n}\left\{ \max_{i<j\leq n,s\in\S}|R_{s}(\bs^{i},\bs^{j})|\leq\delta\right\} >e^{-t(n-1)N}.
\]
Hence, $q\equiv0$ is multi-samplable.\qed

\subsection{Proof of Lemma \ref{lem:transfertoSN}}

We may identify
\[
B(m,0)=\left\{ \bs\in S_{N}:\,\forall s\in\S,\,R_{s}(\bs,m)=R_{s}(m,m)\right\} 
\]
with the product of the spheres
\[
\Big\{(\sigma_{i})_{i\in I_{s}}\in S(N_{s}):\,\sum_{i\in I_{s}}\sigma_{i}m_{i}=\sum_{i\in I_{s}}m_{i}^{2}\Big\}.
\]
Endow each of those spheres with the uniform measure, and let $\nu$
be the product measure on $B(m,0)$. Define 
\[
\bar{F}_{N}(m,n,\rho):=\frac{1}{Nn}\log\int_{B(m,n,0,\rho)}e^{\sum_{i=1}^{n}\big(H_{N}(\bs^{i})-H_{N}(m)\big)}d\nu(\bs^{1})\cdots d\nu(\bs^{n}).
\]

We will prove the following two lemmas below. 
\begin{lem}
\label{lem:FbarF}Let $q\in[0,1)^{\S}$ and $m\in S_{N}(q)$ and fix
some $t>0$. Then, for any $\delta,\rho>0$ there exist $\delta_{0}=\delta_{0}(\delta,\rho)>0$
and $\rho_{0}=\rho_{0}(\delta,\rho)>0$ such that if $\rho'<\rho_{0}$
and $\delta'<\delta_{0}$, then for any $n\geq1$ and large $N$,
\begin{align}
\E F_{N}(m,n,\delta,\rho) & \geq\E\bar{F}_{N}(m,n,\rho')+\frac{1}{N}\log\mu(B(m,\delta'))-t,\label{eq:Fbar1}\\
\E F_{N}(m,n,\delta',\rho') & \leq\E\bar{F}_{N}(m,n,\rho)+\frac{1}{N}\log\mu(B(m,\delta'))+t.\label{eq:Fbar2}
\end{align}
\end{lem}

\begin{lem}
\label{lem:FbarFtilde}Let $q\in[0,1)^{\S}$ and $m\in S_{N}(q)$.
Then for $n\geq1$ and $\rho>0$,
\[
\lim_{N\to\infty}\big|\E\bar{F}_{N}(m,n,\rho)-\E\tilde{F}_{N}^{q}(n,\rho)\big|=0.
\]
\end{lem}

Let $t>0$ and let $n$ be some large number and $\delta$ and $\rho$
be some small numbers. By Lemma \ref{lem:FbarF} there exist $\delta'<\delta$
and $\rho'<\rho$ such that for large $N$,
\begin{align*}
\E F_{N}(m,n,\delta,\rho) & \geq\E\bar{F}_{N}(m,n,\rho')+\frac{1}{N}\log\mu(B(m,\delta'))-t.
\end{align*}
 Assuming $\delta$ is small enough,
\[
\limsup_{N\to\infty}\Big|\frac{1}{N}\log\mu(B(m,\delta'))-\frac{1}{2}\sum_{s\in\S}\lambda_{s}\log(1-q(s))\Big|<t.
\]
Combining this with Lemma \ref{lem:FbarFtilde} we have that for large
$N$,
\begin{align*}
\E F_{N}(m,n,\delta,\rho) & \geq\E\tilde{F}_{N}^{q}(n,\rho')+\frac{1}{2}\sum_{s\in\S}\lambda_{s}\log(1-q(s))-3t.
\end{align*}

Note that by Lemmas \ref{lem:remove1spin} and \ref{lem:manyto1replica},
assuming $\rho'<\rho$ are small enough,
\[
\limsup_{N\to\infty}\Big|\E\tilde{F}_{N}^{q}(n,\rho')-\E\tilde{F}_{N}^{q}(n,\rho)\Big|<t.
\]
Hence,
\begin{align*}
\E F_{N}(m,n,\delta,\rho) & \geq\E\tilde{F}_{N}^{q}(n,\rho)+\frac{1}{2}\sum_{s\in\S}\lambda_{s}\log(1-q(s))-4t.
\end{align*}

By a similar argument, for large $n$ and small $\delta$ and $\rho$,
\begin{align*}
\E F_{N}(m,n,\delta,\rho) & \leq\E\tilde{F}_{N}^{q}(n,\rho)+\frac{1}{2}\sum_{s\in\S}\lambda_{s}\log(1-q(s))+4t.
\end{align*}
This completes the proof of Lemma \ref{lem:transfertoSN}. It remains
to prove Lemmas \ref{lem:FbarF} and \ref{lem:FbarFtilde}.

\subsubsection{Proof of Lemma \ref{lem:FbarFtilde}}

The lemma is an easy consequence of the fact that the mapping $\bs\mapsto\tilde{\bs}$
as defined in (\ref{eq:sigmatilde}) maps $B(m,0)$ bijectively to
\begin{equation}
\Big\{\bs\in S_{N}:\,s\in\S,\,R_{s}(\bs,m)=0\Big\}=\prod_{s\in\S}\Big\{(\sigma_{i})_{i\in I_{s}}:\,\sum_{i\in I_{s}}\sigma_{i}^{2}=N_{s},\sum_{i\in I_{s}}\sigma_{i}m_{i}=0\Big\},\label{eq:set}
\end{equation}
and satisfies (\ref{eq:xiqcov}). \qed

\subsubsection{Proof of Lemma \ref{lem:FbarF}}

Fix $q\in(0,1)^{\S}$ and $m\in S_{N}(q)$, where, for simplicity,
we assume for now that $q(s)>0$ for all $s\in\S$. Define the mapping
\[
\bs\mapsto\varphi(\bs)=\varphi(\bs,m)\in B(m,0)
\]
by first defining the projection $\btau=(\tau_{i})_{i\leq N}$ by
\[
\tau_{i}:=\sigma_{i}-\frac{R_{s}(\bs,m)}{R_{s}(m,m)}m_{i},\quad\text{if }i\in I_{s},
\]
and then defining $\varphi(\bs)=\bp=(\pi_{i})_{i\leq N}$ by
\begin{equation}
\pi_{i}:=m_{i}+\sqrt{\frac{1-R_{s}(m,m)}{R_{s}(\btau,\btau)}}\tau_{i},\quad\text{if }i\in I_{s}.\label{eq:pii}
\end{equation}
We assume here that $\bs$ is such $R_{s}(\btau,\btau)>0$. Any $\bs\in B(m,\delta)$
satisfies this if we assume, as we will from now on, that $\delta<\min_{s}(\sqrt{q(s)}-q(s))$.

Define the free energy
\[
\hat{F}_{N}(m,n,\delta,\rho)=\frac{1}{Nn}\log\int_{\hat{B}(m,n,\delta,\rho)}e^{\sum_{i=1}^{n}\big(H_{N}(\bs^{i})-H_{N}(m)\big)}d\mu(\bs^{1})\cdots d\mu(\bs^{n}),
\]
where
\[
\hat{B}(m,n,\delta,\rho)=\left\{ (\bs^{i})_{i\leq n}\in B(m,\delta)^{n}:\,\forall i\neq j,\,s\in\S,\,\big|R_{s}(\varphi(\bs^{i}),\varphi(\bs^{j}))-R_{s}(m,m)\big|\leq\rho\right\} .
\]

One can check that 
\begin{equation}
\bs\in B(m,\delta)\Longrightarrow R_{s}(\varphi(\bs)-\bs,\varphi(\bs)-\bs)\leq\frac{2\delta}{\sqrt{q(s)}},\label{eq:phidist}
\end{equation}
and therefore
\begin{equation}
\bs^{1},\bs^{2}\in B(m,\delta)\Longrightarrow|R_{s}(\varphi(\bs^{1}),\varphi(\bs^{2}))-R_{s}(\bs^{1},\bs^{2})|\leq\frac{\sqrt{8\delta}}{q(s)^{1/4}}.\label{eq:twoptdist}
\end{equation}
Thus, for $\alpha:=\max_{s\in\S}q(s)^{-1/4}$,
\begin{align*}
B(m,n,\delta,\rho) & \subset\hat{B}(m,n,\delta,\rho+\alpha\sqrt{8\delta}),\\
\hat{B}(m,n,\delta,\rho) & \subset B(m,n,\delta,\rho+\alpha\sqrt{8\delta}),
\end{align*}
and
\begin{align*}
F_{N}(m,n,\delta,\rho) & \leq\hat{F}_{N}(m,n,\delta,\rho+\alpha\sqrt{8\delta}),\\
\hat{F}_{N}(m,n,\delta,\rho) & \leq F_{N}(m,n,\delta,\rho+\alpha\sqrt{8\delta}).
\end{align*}

For any $t=(t_{s})\in(-\delta,\delta)^{\S}$, define the vector $m(t)=(m_{i}(t))_{i\leq N}$
by 
\[
m_{i}(t):=\Big(1+\frac{t_{s}}{R_{s}(m,m)}\Big)m_{i},\quad\text{if }i\in I_{s},
\]
and define $\varphi_{t}(\bs)=\varphi(\bs,m(t))\in B(m(t),0)$. Note
that 
\[
\hat{B}(m,n,\delta,\rho)=\Big\{(\varphi_{t^{i}}(\bs^{i}))_{i\leq n}:\,(\bs^{i})_{i\leq n}\in B(m,n,0,\rho),\,t^{1},\ldots,t^{n}\in(-\delta,\delta)^{\S}\Big\}.
\]
By the co-area formula, we may write
\begin{align*}
 & \hat{F}_{N}(m,n,\delta,\rho)\\
 & =\frac{1}{Nn}\log\int_{B(m,n,0,\rho)}e^{\sum_{i=1}^{n}\big(H_{N}(\bs^{i})-H_{N}(m)\big)}\Theta(\bs^{1},\ldots,\bs^{n})d\nu(\bs^{1})\cdots d\nu(\bs^{n}),
\end{align*}
where
\[
\Theta(\bs^{1},\ldots,\bs^{n})=\int_{(-\delta,\delta)^{\S}}\cdots\int_{(-\delta,\delta)^{\S}}e^{\sum_{i=1}^{n}\big(H_{N}(\varphi_{t^{i}}(\bs^{i}))-H_{N}(\bs^{i})\big)}D(t^{1},\ldots,t^{n})dt^{1}\cdots dt^{n}
\]
and $D(t^{1},\ldots,t^{n})$ is a Jacobian which integrates to 
\[
\int_{(-\delta,\delta)^{\S}}\cdots\int_{(-\delta,\delta)^{\S}}D(t^{1},\ldots,t^{n})dt^{1}\cdots dt^{n}=\frac{\mu^{\otimes n}(\hat{B}(m,n,\delta,\rho))}{\nu^{\otimes n}(B(m,n,0,\rho))}=\mu(B(m,\delta))^{n}.
\]
For $\bs^{i}\in B(m,\delta)$ and $t^{i}\in(-\delta,\delta)^{\S}$,
$\varphi_{t^{i}}(\bs^{i})\in B(m,\delta)$ and $\varphi(\varphi_{t^{i}}(\bs^{i}))=\varphi(\bs^{i})$.
By (\ref{eq:phidist}), 
\begin{align*}
 & R_{s}(\varphi_{t^{i}}(\bs^{i})-\bs^{i},\varphi_{t^{i}}(\bs^{i})-\bs^{i})\\
 & \leq\Big(R_{s}(\varphi_{t^{i}}(\bs^{i})-\varphi(\bs^{i}),\varphi_{t^{i}}(\bs^{i})-\varphi(\bs^{i}))^{1/2}\\
 & +R_{s}(\varphi(\bs^{i})-\bs^{i},\varphi(\bs^{i})-\bs^{i})^{1/2}\Big)^{2}\leq\frac{8\delta}{\sqrt{q(s)}}.
\end{align*}

Hence, if the event in (\ref{eq:LipProb}) occurs for some $L$, then
for any $(\bs^{1},\ldots,\bs^{n})\in B(m,n,0,\rho)$,
\[
\Big|\sum_{i=1}^{n}\big(H_{N}(\varphi_{t^{i}}(\bs^{i}))-H_{N}(\bs^{i})\big)\Big|\leq nL\max_{s\in\S}\sqrt{\frac{8\delta}{\sqrt{q(s)}}}=n\sqrt{8\delta}\alpha L
\]
and
\[
\Big|\bar{F}_{N}(m,n,\rho)+\frac{1}{N}\log\mu(B(m,\delta))-\hat{F}_{N}(m,n,\delta,\rho)\Big|\leq\sqrt{8\delta}\alpha L.
\]

Given $t$, $\delta$ and $\rho$, we may choose some smaller $\delta'$
and $\rho'$ such that $\rho>\rho'+\alpha\sqrt{8\delta'}$ and therefore
from the above,
\begin{align*}
F_{N}(m,n,\delta,\rho) & \geq F_{N}(m,n,\delta',\rho'+\alpha\sqrt{8\delta'})\geq\hat{F}_{N}(m,n,\delta',\rho')\\
 & \geq\bar{F}_{N}(m,n,\rho')+\frac{1}{N}\log\mu(B(m,\delta'))-\sqrt{8\delta'}\alpha L.
\end{align*}
Assume that $L$ is large enough so that by Lemma \ref{lem:Lipschitz}
the probability in (\ref{eq:LipProb}) goes to $1$ with $N$ and
that $\delta'$ is small enough so that $\sqrt{8\delta'}\alpha L<t/2$.
Then, (\ref{eq:Fbar1}) follows from the concentration of the free
energies $F_{N}(m,n,\delta,\rho)$ and $\bar{F}_{N}(m,n,\rho')$ around
their mean (e.g., by \cite[Theorem 1.2]{PanchenkoBook}). A similar
argument gives (\ref{eq:Fbar2}). 

Lastly, recall that we assumed that $q(s)>0$. If $q(s)=0$ for some
$s\in\S$, then one only needs to redefine $\varphi(\bs)=\bp$ by
setting $\pi_{i}=\sigma_{i}$ instead of (\ref{eq:pii}) for all $i\in I_{s}$
and such $s$. The inequalities (\ref{eq:phidist}) and (\ref{eq:twoptdist})
then become trivial for such $s$, and up to obvious modifications
the proof remains the same.\qed

\section{\label{sec:pfMAin}TAP representation: Proof of Theorems \ref{thm:TAP}
and \ref{thm:Onsager}}

Theorem \ref{thm:TAP} follows directly from Propositions \ref{prop:TAPwithoutcorrection}
and \ref{prop:correction}.  Theorem \ref{thm:Onsager} follows from
the following two results.
\begin{cor}
\label{cor:qRS}Let $q\in[0,1)^{\S}$ be a maximal multi-samplable
overlap of the mixture $\xi$, and assume $\xi$ satisfies the convergence
condition (\ref{eq:conv}). Then, for any $s\in\S$ and $\tau>0$,
\begin{equation}
\limsup_{N\to\infty}\frac{1}{N}\log\E\max_{\bp\in S_{N}}G_{N}^{q}\left\{ |R_{s}(\bs,\bp)|\geq\tau\right\} <0,\label{eq:RS}
\end{equation}
where $G_{N}^{q}$ is the Gibbs measure corresponding to the Hamiltonian
$H_{N}^{q}(\bs)$ with mixture $\xi_{q}(x)$.
\end{cor}

\begin{lem}
\label{lem:RS}Let $\xi$ be some mixture. If for any $s\in\S$ and
$\tau>0$, 
\begin{equation}
\limsup_{N\to\infty}\frac{1}{N}\log\E\max_{\bp\in S_{N}}G_{N}\left\{ |R_{s}(\bs,\bp)|\geq\tau\right\} <0,\label{eq:temp1}
\end{equation}
then 
\[
\E F_{N}=\frac{1}{2}\xi(1)+o(1).
\]
\end{lem}

In the rest of the section we will prove the results above, and the
following two results. 
\begin{lem}
\label{lem:nesting}If $q$ is a multi-samplable overlap of $\xi$
and $q'$ is a multi-samplable overlap of $\xi_{q}$, then
\[
q+(1-q)q'
\]
is also a multi-samplable overlap of $\xi$. Therefore, if $q\in[0,1)^{\S}$
is a maximal multi-samplable overlap for some mixture $\xi$, then
the only multi-samplable overlap of $\xi_{q}$ is $q'\equiv0$.
\end{lem}

\begin{lem}
\label{lem:q0RS}If for some $s\in\S$ and $\tau>0$,
\begin{equation}
\limsup_{N\to\infty}\frac{1}{N}\log\E\max_{\bp\in S_{N}}G_{N}\left\{ R_{s}(\bs,\bp)\geq\tau\right\} =0,\label{eq:cont}
\end{equation}
then for some $q\in[0,1)^{\S}$ such that $q(s)\geq\tau^{2}/2$ and
some subsequence $N_{\ell}$, as $\ell\to\infty$,
\begin{equation}
\E F_{N_{\ell}}=\E\EsNell(q)+\frac{1}{2}\sum_{s\in\S}\lambda_{s}\log(1-q(s))+\E F_{N_{\ell}}(q)+o(1).\label{eq:TAPsubseq}
\end{equation}
\end{lem}

\subsection{Proof of Lemma \ref{lem:nesting}}

Let $q$ and $q'$ be multi-samplable overlaps as in the lemma, and
define $\hat{q}(s)=q(s)+(1-q(s))q'(s)$. By Theorem \ref{thm:TAP},
\begin{equation}
\begin{aligned}\E F_{N} & =\E\EsN(q)+\frac{1}{2}\sum_{s\in\S}\log(1-q(s))+\E F_{N}(q)+o(1),\\
\E F_{N}(q) & =\E\EsN^{q}(q')+\frac{1}{2}\sum_{s\in\S}\log(1-q'(s))+\E F_{N}^{q}(q')+o(1),
\end{aligned}
\label{eq:TAPqq'}
\end{equation}
where the terms in the first line are defined using $\xi$, and terms
in the second line, with superscript $q$, are defined similarly using
the mixture $\xi_{q}$. Moreover, by the same theorem, to prove the
lemma we need to show that
\begin{align}
\E F_{N} & \leq\E\EsN(\hat{q})+\frac{1}{2}\sum_{s\in\S}\log(1-\hat{q}(s))+\E F_{N}(\hat{q})+o(1).\label{eq:TAPqhat}
\end{align}

Since
\[
1-\hat{q}(s)=(1-q(s))(1-q'(s)),
\]
the sum of the two logarithmic terms in (\ref{eq:TAPqq'}) is equal
to the logarithmic term in (\ref{eq:TAPqhat}). Also, from the definition
(\ref{eq:xiq}) of $\xi_{q}$, it is straightforward to verify that
$F(\hat{q})$ and $F^{q}(q')$ are the free energies of the same mixture
$\xi_{\hat{q}}$ and are therefore equal. Hence, to prove (\ref{eq:TAPqhat})
it remains to show that for any $q$, $q'$ and $\hat{q}$ defined
as above, 
\begin{equation}
\E\EsN(q)+\E\EsN^{q}(q')\leq\E\EsN(\hat{q})+o(1).\label{eq:Esineq}
\end{equation}

This follows by a straightforward modification of the argument used
for the single-species models in the proof of (6.4) in Lemma 31 of
\cite{FElandscape}. Here, in the multi-species case, one needs to
use the decomposition as in (\ref{eq:1spindecom}) and the Lipschitz
property of Lemma \ref{lem:Lipschitz} which generalized the analogous
results for the single-species used in \cite{FElandscape}. (We remark
that (\ref{eq:Esineq}) holds for any $q$ and $q'$ in $[0,1)^{\S}$,
not necessarily multi-samplable overlaps.)\qed

\subsection{Proof of Lemma \ref{lem:q0RS}}

Let $\tau>0$ and $s'\in\S$. Let $n\geq1$ and $\epsilon>0$ such
that $1/\epsilon$ is an integer and $\frac{1}{n}+\epsilon\leq\tau^{2}/2$.
Denote 
\[
Q(\epsilon)=\left\{ q:\,\forall s,\,q(s)\in\{0,\epsilon,\ldots,1-\epsilon\}\right\} .
\]
Decompose $[0,1)^{\S}$ into $\epsilon^{-|\S|}$ disjoint sets
\begin{equation}
[0,1)^{\S}=\bigcup_{q\in Q(\epsilon)}\prod_{s\in\S}[q(s),q(s)+\epsilon).\label{eq:dec}
\end{equation}

By Ramsey's theorem, if $k=k(n)$ is large enough, then for any choice
of $\bs^{1},\ldots,\bs^{k}\in S_{N}$, there is a subset of size $n$
such that the overlaps $R(\bs^{i},\bs^{j})$ of all pairs from it
belong to the same subset from the decomposition (\ref{eq:dec}),
and therefore for some $q\in Q(\epsilon)$,
\begin{equation}
\max_{s\in\S}|R_{s}(\bs^{i},\bs^{j})-q(s)|<\epsilon.\label{eq:Rsq}
\end{equation}

Suppose that $\bs^{1},\ldots,\bs^{n}\in S_{N}$ and $q\in Q(\epsilon)$
are as above, i.e. for any $i\neq j$, (\ref{eq:Rsq}) holds, and
that for some $\bp\in S_{N}$ and every $i\leq n$, $R_{s'}(\bs^{i},\bp)>\tau$.
Then, for $m=\frac{1}{n}\sum_{i\leq n}\bs^{i}$,
\begin{align*}
\Big|R_{s'}(m,m)-\Big(\frac{1}{n}+\frac{n-1}{n}q(s')\Big)\Big| & <\epsilon,
\end{align*}
and by Cauchy\textendash Schwarz,
\[
R_{s'}(m,\bp)>\tau\implies R_{s'}(m,m)>\tau^{2}.
\]
Hence, since we assumed that $\frac{1}{n}+\epsilon\leq\tau^{2}/2$,
\begin{equation}
q(s')\geq\tau^{2}/2.\label{eq:qs'}
\end{equation}

We therefore have that
\begin{equation}
\indic\{(\bs^{i_{1}},\ldots,\bs^{k})\in A_{k}(\bp,\tau)\}\leq\sum\indic\{(\bs^{i_{1}},\ldots,\bs^{i_{n}})\in B_{n}(q,\epsilon)\}\label{eq:nk}
\end{equation}
where the sum is over all $q\in Q(\epsilon)$ with $q(s')\geq\tau^{2}/2$
and $1\leq i_{1}<\cdots<i_{n}\leq k$ and
\begin{align*}
A_{k}(\bp,\tau) & =\left\{ (\bs^{i})_{i\leq k}\in S_{N}^{k}:\,R_{s'}(\bs^{i},\bp)>\tau\right\} ,\\
B_{n}(q,\epsilon) & =\left\{ (\bs^{i})_{i\leq n}\in S_{N}^{n}:\,\max_{s\in\S}|R_{s}(\bs^{i},\bs^{j})-q(s)|<\epsilon,\,\forall i\neq j\right\} .
\end{align*}
By taking the expectation w.r.t. $G_{N}^{\otimes k}$ of both sides
of (\ref{eq:nk}),
\begin{equation}
\begin{aligned} & \binom{k}{n}\epsilon^{-|\S|}\max_{q\in Q(\epsilon):\,q(s')>\frac{\tau^{2}}{2}}G_{N}^{\otimes n}\Big\{\max_{s\in\S,i<j\leq n}|R_{s}(\bs^{i},\bs^{j})-q(s)|<\epsilon\Big\}\\
 & \geq\max_{\bp\in S_{N}}G_{N}^{\otimes k}\Big\{\forall i\leq k,\,R_{s'}(\bs^{i},\bp)>\tau\Big\}=\Big(\max_{\bp\in S_{N}}G_{N}\left\{ R_{s'}(\bs,\bp)>\tau\right\} \Big)^{k}.
\end{aligned}
\label{eq:posgoodq}
\end{equation}

Now, assume that
\[
\limsup_{N\to\infty}\frac{1}{N}\log\E\max_{\bp\in S_{N}}G_{N}\left\{ R_{s'}(\bs,\bp)\geq\tau\right\} =0.
\]
Since $Q(\epsilon)$ is finite, from (\ref{eq:posgoodq}) we obtain
that there exists $q=q(n,\epsilon)\in Q(\epsilon)$ independent of
$N$, such that $q(s')\geq\tau^{2}/2$ and
\[
\limsup_{N\to\infty}\frac{1}{N}\log\E G_{N}^{\otimes n}\left\{ \max_{s\in\S,i<j\leq n}|R_{s}(\bs^{i},\bs^{j})-q(s)|<\epsilon\right\} =0.
\]
If $q_{\star}$ is a limit point of $q(n_{i},\epsilon_{i})$ for some
sequences $n_{i}\to\infty$ and $\epsilon_{i}\to0$, then $q_{\star}(s')\geq\tau^{2}/2$
and for all $n$ and $\epsilon$,
\[
\limsup_{N\to\infty}\frac{1}{N}\log\E G_{N}^{\otimes n}\left\{ \max_{s\in\S,i<j\leq n}|R_{s}(\bs^{i},\bs^{j})-q_{\star}(s)|<\epsilon\right\} =0.
\]
Since the expectation above is increasing in $\epsilon$ and decreasing
in $n$, by diagonalization there exists a subsequence $N_{\ell}$
such that for any $n$ and $\epsilon$,
\[
\lim_{\ell\to\infty}\frac{1}{N_{\ell}}\log\E G_{N_{\ell}}^{\otimes n}\left\{ \max_{s\in\S,i<j\leq n}|R_{s}(\bs^{i},\bs^{j})-q_{\star}(s)|<\epsilon\right\} =0.
\]

By going over its proof, one can check that Theorem \ref{thm:TAP}
also holds for subsequences in the sense that an overlap $q\in[0,1)^{\S}$
satisfies (\ref{eq:multisamp-1}) for any $n$ and $\epsilon$ on
a subsequence $N_{\ell}$ if and only if (\ref{eq:TAPrep}) holds
on the same subsequence. Thus, $q_{\star}$ satisfies (\ref{eq:TAPsubseq}).\qed

\subsection{Proof of Corollary \ref{cor:qRS}}

Suppose that $q\in[0,1)^{\S}$ is a maximal multi-samplable overlap
of a mixture $\xi$ which satisfies the convergence condition (\ref{eq:conv}).
Assume towards contradiction that for some $s\in\S$ and $\tau>0$,
\begin{equation}
\limsup_{N\to\infty}\frac{1}{N}\log\E\max_{\bp\in S_{N}}G_{N}^{q}\left\{ |R_{s}(\bs,\bp)|\geq\tau\right\} =0.\label{eq:cont-1}
\end{equation}

Since
\[
G_{N}^{q}\left\{ |R_{s}(\bs,\bp)|\geq\tau\right\} =G_{N}^{q}\left\{ R_{s}(\bs,\bp)\geq\tau\right\} +G_{N}^{q}\left\{ R_{s}(\bs,-\bp)\geq\tau\right\} ,
\]
we also have that
\[
\limsup_{N\to\infty}\frac{1}{N}\log\E\max_{\bp\in S_{N}}G_{N}^{q}\left\{ R_{s}(\bs,\bp)\geq\tau\right\} =0.
\]

By Lemma \ref{lem:q0RS}, for some $q'$ with $q'(s)\geq\tau^{2}/2$
and $N_{\ell}$ as in the lemma, as $\ell\to\infty$,
\[
\E F_{N_{\ell}}(q)=\E F_{N_{\ell}}^{q}=\E\EsNell^{q}(q')+\frac{1}{2}\sum_{s\in\S}\lambda_{s}\log(1-q'(s))+\E F_{N_{\ell}}^{q}(q')+o(1),
\]
where all the quantities above with superscript $q$ are defined using
the Hamiltonian $H_{N}^{q}(\bs)$ with mixture $\xi_{q}(x)$.

Since $q$ is multi-samplable, by the same argument we used in the
proof of Lemma \ref{lem:nesting}, defining $\hat{q}(s)=q(s)+(1-q(s))q'(s)$,
on the same subsequence $N_{\ell}$,
\[
\E F_{N_{\ell}}=\E\EsNell(\hat{q})+\frac{1}{2}\sum_{s\in\S}\log(1-\hat{q}(s))+\E F_{N_{\ell}}(\hat{q})+o(1).
\]
From (\ref{eq:conv}), we have that the same holds without moving
to a subsequence. That is,
\[
\E F_{N}=\E\EsN(\hat{q})+\frac{1}{2}\sum_{s\in\S}\log(1-\hat{q}(s))+\E F_{N}(\hat{q})+o(1).
\]

By Theorem \ref{thm:TAP}, $\hat{q}$ is therefore multi-samplable.
This contradicts the maximality of $q$, and we conclude that (\ref{eq:cont-1})
does not hold for any $s\in\S$ and $\tau>0$.\qed

\subsection{Proof of Lemma \ref{lem:RS}}

From Jensen's inequality, $\E F_{N}\leq\frac{1}{2}\xi(1)$. Hence,
to prove the lemma it will be enough to show that 
\[
\liminf_{N\to\infty}\E\log Z_{N+1}-\E\log Z_{N}\geq\frac{1}{2}\xi(1).
\]

We remark that changing the values of $N_{s}$ while keeping the same
limiting proportions $\lambda_{s}$ results in changing the free energy
$\E F_{N}$ by an amount which vanishes as $N\to\infty$.\footnote{\label{fn:Ns} Note that given $(N_{s})$ and $(N'_{s})$ with the
same $N\to\infty$ limit, one can relate each of the associated free
energies to that of the model corresponding to $\min\{N_{s},N_{s}'\}$,
by integrating over the `extra' coordinates. Using the Lipschitz property
in Lemma \ref{lem:Lipschitz} and the concentration of the free energy,
one can then see that the difference of the free energies goes to
$0$ as $N\to\infty$. } By a similar argument as in the footnote, one can also verify that
(\ref{eq:temp1}) also holds if we change the values of $N_{s}$,
as long as the limits $\lambda_{s}$ are the same. Therefore we may,
and will, assume that the values $N_{s}$ are non-decreasing with
$N$. In this case, as $N$ increases to $N+1$, $N_{s}$ increases
to $N_{s}+1$ for exactly one $s\in\S$ while all other $N_{s}$ remain
the same. We will denote the species $s\in\S$ that increases by $s_{\star}=s_{\star}(N)$.

Instead of working with $S_{N}$ in this proof we will work with 
\[
T_{N}=\left\{ (\sigma_{1},\ldots,\sigma_{N})\in\R^{N}:\,\forall s\in\S,\,\sum_{i\in I_{s}}\sigma_{i}^{2}=1\right\} ,
\]
by defining for $\bs\in T_{N}$, 
\begin{equation}
h_{N}(\bs)=H_{N}(\tilde{\bs})=\sqrt{N}\sum_{k=1}^{\infty}\sum_{i_{1},\dots,i_{k}=1}^{N}\bar{\Delta}_{i_{1},\dots,i_{k}}J_{i_{1},\dots,i_{k}}\sigma_{i_{1}}\cdots\sigma_{i_{k}}\label{eq:hN}
\end{equation}
where for any $s\in\S$ and $i\in I_{s}$,$\tilde{\sigma}_{i}=\sqrt{N_{s}}\sigma_{i}$
and if $\#\{j\leq k:\,i_{j}\in I_{s}\}=p(s)$ for any $s\in\S$, then
$\Delta_{i_{1},\dots,i_{k}}=\Delta_{i_{1},\dots,i_{k}}(N)$ is defined
by
\[
\bar{\Delta}_{i_{1},\dots,i_{k}}^{2}=\Delta_{p}^{2}\frac{\prod_{s\in\S}p(s)!}{k!}\prod_{s\in\S}.
\]

Let $N\geq1$. When the system size is $N$, assume WLOG that $I_{s_{\star}}=\{N-N_{s_{\star}}+1,\ldots,N\}$
and when the system size is $N+1$ assume that $I_{s_{\star}}=\{N-N_{s_{\star}}+1,\ldots,N+1\}$.
We will use $\br$ to denote points from $T_{N+1}$ and $\bs$ for
points from $T_{N}$, and will assume everywhere that $\br$, $\bs$
and $\varepsilon$ are related to each other by 
\begin{equation}
\rho_{i}=\begin{cases}
\sigma_{i} & \text{if }i\leq N-N_{s_{\star}},\\
\sqrt{1-\varepsilon^{2}}\sigma_{i} & \text{if }N-N_{s_{\star}}<i\leq N,\\
\varepsilon & \text{if }i=N+1.
\end{cases}\label{eq:rhobs}
\end{equation}
Here $\varepsilon$ plays the role of the usual cavity coordinate.

Using the same variables $J_{i_{1},\ldots,i_{k}}$ as in (\ref{eq:hN}),
define 
\[
h_{N+1}^{(1)}(\br)=\sqrt{N}\sum_{k=1}^{\infty}\sum_{i_{1},\dots,i_{k}=1}^{N}\bar{\Delta}_{i_{1},\dots,i_{k}}J_{i_{1},\dots,i_{k}}\sigma_{i_{1}}\cdots\sigma_{i_{k}}(1-\varepsilon^{2})^{\frac{1}{2}\#\{j\leq k:\,i_{j}\in I_{s_{\star}}\}}.
\]
If $\bx=(x_{i})_{i\leq N}$ is defined by $x_{i}=\rho_{i}$, then
\[
\E h_{N+1}^{(1)}(\br)h_{N+1}^{(1)}(\br')=N\xi(Q(\bx,\bx')),
\]
where we define $(Q(\bx,\bx'))(s):=\sum_{i\in I_{s}}x_{i}x_{i}'$.

If we write $h_{N+1}(\br)$ as a polynomial in $\br$ using (\ref{eq:hN}),
then $\sqrt{\frac{N+1}{N}}h_{N+1}^{(1)}(\br)$ is the sum of all terms
which do not contain the last coordinate $\varepsilon$. Hence, we
may assume, as we will, that $h_{N+1}(\br)$, $h_{N+1}^{(1)}(\br)$
and an additional process $h_{N+1}^{(2)}(\br)$ independent of $h_{N+1}^{(1)}(\br)$
are defined on the same probability space, such that 
\[
h_{N+1}(\br)=h_{N+1}^{(1)}(\br)+h_{N+1}^{(2)}(\br),
\]
and 
\[
\begin{aligned} & \E\Big\{ h_{N+1}^{(2)}(\br)h_{N+1}^{(2)}(\br')\Big\}\\
 & =\E\Big\{ h_{N+1}(\br)h_{N+1}(\br')\Big\}-\E\Big\{ h_{N+1}^{(1)}(\br)h_{N+1}^{(1)}(\br')\Big\}.
\end{aligned}
\]

By a Taylor approximation one can check that, if $\varepsilon,\,\varepsilon'\in(-t,t)$,
\begin{equation}
\begin{aligned} & \Big|\E\Big\{ h_{N+1}^{(2)}(\br)h_{N+1}^{(2)}(\br')\Big\}-\xi(Q(\bs,\bs'))-N\varepsilon\varepsilon'\frac{d}{dx(s_{\star})}\xi(Q(\bs,\bs'))\Big|\\
 & \leq4Nt^{4}\frac{d^{2}}{dx(s_{\star})^{2}}\xi(1)+4t^{2}\frac{d}{dx(s_{\star})}\xi(1).
\end{aligned}
\label{eq:h2cov}
\end{equation}

Fix $C>0$ and define $B_{N}:=\{\br\in T_{N+1}:\,|\varepsilon|<N^{-1/2}C\}$.
Note that
\begin{align*}
\log Z_{N+1}-\log Z_{N} & =\log\int_{T_{N+1}}e^{h_{N+1}(\br)}d\mu(\br)-\log\int_{T_{N}}e^{h_{N}(\bs)}d\mu(\bs)\\
 & \geq\log\int_{B_{N}}e^{h_{N+1}(\br)}d\mu(\br)-\log\int_{T_{N}}e^{h_{N}(\bs)}d\mu(\bs)=:W_{N},
\end{align*}
where $\mu=\mu_{N}$ is the product of the uniform measures on the
spheres corresponding to the subsets of coordinates $I_{s}$. We will
show that for any $t>0$, if $C$ is large enough, 
\begin{equation}
\lim_{N\to\infty}\P\left(W_{N}>\frac{1}{2}\xi(1)-t\right)=1.\label{eq:WNbound}
\end{equation}

Recall (\ref{eq:h2cov}) and our assumption that (\ref{eq:RS}) holds.
Combining the latter with a straightforward modification of the argument
used to prove Lemma 34 in \cite{FElandscape} the following lemma
follows. 
\begin{lem}
\label{lem:h(2)}For any $\delta>0$, with probability going to $1$
as $N\to\infty$,
\begin{align*}
\int_{B_{N}}e^{h_{N+1}(\br)}d\br & \geq(1-\delta/2)\int_{B_{N}}e^{h_{N+1}^{(1)}(\br)+\frac{1}{2}\text{Var}(h_{N+1}^{(2)}(\br))}d\br\\
 & \geq(1-\delta)\int_{B_{N}}e^{h_{N+1}^{(1)}(\br)+\frac{1}{2}(\xi(1)+N\varepsilon^{2}\frac{d}{dx(s_{\star})}\xi(1))}d\br.
\end{align*}
\end{lem}

Define
\begin{align*}
g_{N}(\bs) & =\sqrt{N}\sum_{k=1}^{\infty}\sum_{i_{1},\dots,i_{k}=1}^{N}p_{i_{1},\dots,i_{k}}(s_{\star})\bar{\Delta}_{i_{1},\dots,i_{k}}J_{i_{1},\dots,i_{k}}\sigma_{i_{1}}\cdots\sigma_{i_{k}},\\
\zeta_{N}(\br) & =\sqrt{N}\sum_{k=1}^{\infty}\sum_{i_{1},\dots,i_{k}=1}^{N}\tau_{i_{1},\dots,i_{k}}(\varepsilon)\bar{\Delta}_{i_{1},\dots,i_{k}}J_{i_{1},\dots,i_{k}}\sigma_{i_{1}}\cdots\sigma_{i_{k}},
\end{align*}
where $p_{i_{1},\dots,i_{k}}(s_{\star}):=\#\{j\leq k:\,i_{j}\in I_{s_{\star}}\}$
and 
\[
\tau_{i_{1},\dots,i_{k}}(\varepsilon):=(1-\varepsilon^{2})^{\frac{1}{2}p_{i_{1},\dots,i_{k}}(s_{\star})}-1+\frac{1}{2}p_{i_{1},\dots,i_{k}}(s_{\star})\varepsilon^{2}.
\]
Then, 
\begin{equation}
h_{N+1}^{(1)}(\br)=h_{N}(\bs)-\frac{1}{2}\varepsilon^{2}g_{N}(\bs)+\zeta_{N}(\br).\label{eq:hN+1}
\end{equation}

Note that 
\[
\frac{1}{2}\frac{d}{dt}\Big|_{t=0}\E\Big((h_{N}(\bs)+tg_{N}(\bs))^{2}\Big)=N\frac{d}{dx(s_{\star})}\xi(1).
\]
Therefore, by the same argument as in the proof of Lemma 35 in \cite{FElandscape},
we have the following lemma.
\begin{lem}
\label{lem:gN}For any $\delta>0$, with probability going to $1$
as $N\to\infty$, there exists a (random) subset $D_{N}\subset T_{N}$
such that
\[
\int_{D_{N}}e^{h_{N}(\bs)}d\bs\geq(1-\delta)\int_{T_{N}}e^{h_{N}(\bs)}d\bs,
\]
and
\[
\sup_{\bs\in D_{N}}g_{N}(\bs)\leq(1+\delta)N\frac{d}{dx(s_{\star})}\xi(1).
\]
\end{lem}

The last term in (\ref{eq:hN+1}) is negligible by the first bound
in the following lemma.
\begin{lem}
\label{lem:zeta}For some constant $c>0$ depending only on $\xi$
and $C$, for any large enough $t>0$, for large $N$,
\begin{align*}
\P\Big(\sup_{\br\in B_{N}}|\zeta_{N}(\br)| & \geq\frac{t}{N}\Big)<e^{-Nct^{2}},\\
\P\Big(\sup_{\bs\in T_{N}}|g_{N}(\bs)| & \geq Nt\Big)<e^{-Nct^{2}}.
\end{align*}
\end{lem}

\begin{proof}
The bound for the process $g_{N}(\bs)$ follows from Lemma \ref{lem:Lipschitz},
applied with $\bp=0$. For any $i_{1},\dots,i_{k}$ define $\kappa_{i_{1},\dots,i_{k}}:=\#\{j\leq k:\,i_{j}\in I_{s_{\star}}\}$.
For $n\geq1$ define 
\[
\zeta_{N,n}(\bs):=\sqrt{N}\sum_{k=1}^{\infty}\sum_{\substack{i_{1},\dots,i_{k}:\\
p_{i_{1},\dots,i_{k}}(s_{\star})=n
}
}\bar{\Delta}_{i_{1},\dots,i_{k}}J_{i_{1},\dots,i_{k}}\sigma_{i_{1}}\cdots\sigma_{i_{k}}.
\]
By an abuse of notation, write $\zeta_{N}(\bs,\varepsilon)=\zeta_{N}(\br)$
assuming (\ref{eq:rhobs}). Then,
\[
\zeta_{N}(\bs,\varepsilon)=\sum_{k\geq1}\Big((1-\varepsilon^{2})^{\frac{1}{2}k}-1+\frac{1}{2}k\varepsilon^{2}\Big)\zeta_{N,k}(\bs).
\]

By the proof of Lemma \ref{lem:Lipschitz}, specifically (\ref{eq:EsupG}),
\[
\E\sup_{\bs\in K_{N}}\zeta_{N,k}(\bs)=\E\sup_{\bs\in K_{N}}\zeta_{N,k}(\bs)-\zeta_{N,k}(0)\leq2N|\S|\Delta_{k},
\]
where $\Delta_{k}=\sum_{p\in P:\,p(s_{\star})=k}\Delta_{p}^{2}|p|^{4}$.
(Here we work with $\bs\in T_{N}$ instead of $\bs\in S_{N}$ as in
Lemma \ref{lem:Lipschitz}, which in (\ref{eq:EsupG}) amounts to
working with $\Delta_{p}$ corresponding to $\Delta_{i_{1},\dots,i_{k}}$
and not the normalized $\bar{\Delta}_{i_{1},\dots,i_{k}}$.) Since
\[
\Big|(1-\varepsilon^{2})^{\frac{1}{2}k}-1+\frac{1}{2}k\varepsilon^{2}\Big|\leq\frac{1}{4}k(\frac{1}{2}k-1)\varepsilon^{4},
\]
we have that
\[
\E\sup_{\bs\in B_{N}}\zeta_{N}(\bs,\varepsilon)\leq\frac{1}{N}C^{4}|\S|\sum_{k\geq1}k(\frac{1}{2}k-1)\Delta_{k}.
\]

Of course,
\[
\sup_{\bs\in B_{N}}\E\left(\zeta_{N}(\bs,\varepsilon)^{2}\right)\leq C^{8}N^{-3}\sum_{k\geq1}\Big(\frac{1}{4}k(\frac{1}{2}k-1)\Big)^{2}\sum_{p\in P:\,p(s_{\star})=k}\Delta_{p}^{2},
\]
which completes the proof by the Borell-TIS inequality.
\end{proof}
Suppose $f(\br)$ is a smooth function on $T_{N+1}$ and, by an abuse
of notation, write $f(\br)=f(\bs,\varepsilon)$ assuming (\ref{eq:rhobs}).
Then,
\begin{equation}
\begin{aligned}\int_{B_{N}}f(\br)d\br & =\frac{\omega_{N_{s_{\star}}-1}}{\omega_{N_{s_{\star}}}}\int_{T_{N}}\int_{-C/\sqrt{N}}^{C/\sqrt{N}}(1-\varepsilon^{2})^{\frac{N_{s_{\star}}-3}{2}}f(\bs,\varepsilon)d\varepsilon d\mu(\bs)\\
 & =(1+o_{N}(1))\sqrt{\frac{N_{s_{\star}}}{2\pi}}\int_{T_{N}}\int_{-C/\sqrt{N}}^{C/\sqrt{N}}e^{-\frac{N_{s_{\star}}}{2}\varepsilon^{2}}f(\bs,\varepsilon)d\varepsilon d\mu(\bs),
\end{aligned}
\label{eq:coarea}
\end{equation}
where $\omega_{N-1}=2\pi^{N/2}/\Gamma(N/2)$ is the surface area of
$T_{N}$ and $\frac{\omega_{N-1}}{\omega_{N}}=(1+o_{N}(1))\sqrt{\frac{N}{2\pi}}$.
Hence, for any $\delta>0$ and with probability going to $1$,
\begin{align*}
 & \int_{B_{N}}e^{h_{N+1}(\br)}d\br\\
 & \geq(1-2\delta)e^{\frac{1}{2}\xi(1)}\sqrt{\frac{N_{s_{\star}}}{2\pi}}\int_{T_{N}}\int_{-C/\sqrt{N}}^{C/\sqrt{N}}e^{-\frac{N_{s_{\star}}}{2}\varepsilon^{2}+h_{N+1}^{(1)}(\br)+\frac{1}{2}N\varepsilon^{2}\frac{d}{dx(s_{\star})}\xi(1)}d\varepsilon d\mu(\bs)\\
 & \geq(1-3\delta)e^{\frac{1}{2}\xi(1)}\sqrt{\frac{N_{s_{\star}}}{2\pi}}\int_{D_{N}}\int_{-C/\sqrt{N}}^{C/\sqrt{N}}e^{-\frac{N_{s_{\star}}}{2}\varepsilon^{2}+h_{N}(\bs)-\frac{1}{2}\varepsilon^{2}\delta N\frac{d}{dx(s_{\star})}\xi(1)}d\varepsilon d\mu(\bs)\\
 & \geq(1-3\delta)^{2}e^{\frac{1}{2}\xi(1)-\frac{1}{2}C^{2}\delta\frac{d}{dx(s_{\star})}\xi(1)}\P\Big(|X|\leq C\sqrt{N_{s_{\star}}/N}\Big)\int_{T_{N}}e^{h_{N}(\bs)}d\mu(\bs),
\end{align*}
where the first inequality follows from Lemma \ref{lem:h(2)}, the
second and third inequalities follow from (\ref{eq:hN+1}) and Lemmas
\ref{lem:zeta} and \ref{lem:gN}, and $X$ is a standard Gaussian
variable. For any $t$, large enough $C=C(t)$ and small enough $\delta=\delta(t,C)$,
this yields (\ref{eq:WNbound}). The proof of Lemma \ref{lem:RS}
will be completed if we show that for fixed $C$, $W_{N}$ is uniformly
integrable.

From Lemma \ref{lem:zeta} and (\ref{eq:coarea}), since log of the
volume of $B_{N}$ is bounded from below uniformly in $N$, for large
$t>0$,
\begin{equation}
\P\Big(\Big|\log\int_{B_{N}}e^{h_{N+1}^{(1)}(\br)}d\mu(\br)-\log\int_{T_{N}}e^{h_{N}(\bs)}d\mu(\bs)\Big|>a+t\Big)\leq e^{-Nct^{2}}\label{eq:bd1}
\end{equation}
for some constant $a=a(C)>0$. Conditioned on $h_{N+1}^{(1)}(\br)$,
for any $\br$ the variance of $e^{h_{N+1}(\br)}$ is equal to the
unconditional variance of $h_{N+1}^{(2)}(\br)$, which is bounded
by some constant $A=A(C)>0$. Hence, by Jensen's inequality, 
\[
0\leq\E\Big(\log\int_{B_{N}}e^{h_{N+1}(\br)}d\mu(\br)\,\Big|\,h_{N+1}^{(1)}(\br)\Big)-\log\int_{B_{N}}e^{h_{N+1}^{(1)}(\br)}d\mu(\br)\leq\frac{1}{2}A.
\]
Moreover, from the concentration of the free energy (see e.g. \cite[Theorem 1.2]{PanchenkoBook})
applied conditionally on $h_{N+1}^{(1)}(\br)$, 
\begin{equation}
\P\Big(\Big|\log\int_{B_{N}}e^{h_{N+1}(\br)}d\mu(\br)-\log\int_{B_{N}}e^{h_{N+1}^{(1)}(\br)}d\mu(\br)\Big|>\frac{1}{2}A+t\Big)\leq2e^{-\frac{t^{2}}{4A}}.\label{eq:bd2}
\end{equation}
Combining (\ref{eq:bd1}) and (\ref{eq:bd2}) we have that
\[
\P\Big(\Big|W_{N}\Big|>a+\frac{1}{2}A+2t\Big)\leq e^{-Nct^{2}}+2e^{-\frac{t^{2}}{4A}}.
\]
Therefore, $W_{N}$ is uniformly integrable and proof of Lemma \ref{lem:RS}
is completed.\qed

\section*{Appendix: Lipschitz continuity}

Denote the closure of $M_{N}$ by 
\[
\overline{M}_{N}=\left\{ (\sigma_{1},\ldots,\sigma_{N})\in\R^{N}:\,\sum_{i\in I_{s}}\sigma_{i}^{2}\leq N_{s},\,\forall s\in\S\right\} .
\]
The following lemma is proved by an adaptation of the proof of Lemma
6.1 in \cite{TAPChenPanchenkoSubag}. 
\begin{lem}
[Lipschitz continuity] \label{lem:Lipschitz}For any $t>0$ there
exists $L>0$, which depends on $\xi$ and $t$ only, such that for
any $N$,
\begin{equation}
\P\left\{ \forall\bs,\bp\in\overline{M}_{N}:\,\frac{1}{N}|H_{N}(\bs)-H_{N}(\bp)|\leq L\max_{s\in\S}\sqrt{R_{s}(\bs-\bp,\bs-\bp)}\right\} >1-e^{-tN}.\label{eq:LipProb}
\end{equation}
\end{lem}

\begin{proof}
Define 
\begin{align*}
G(u,\bs) & :=u\cdot\nabla H_{N}(\bs)\\
 & :=\sqrt{N}\sum_{k=1}^{\infty}\sum_{i_{1},\dots,i_{k}=1}^{N}\Delta_{i_{1},\dots,i_{k}}J_{i_{1},\dots,i_{k}}(u_{i_{1}}\cdots\sigma_{i_{k}}+\cdots\sigma_{i_{1}}\cdots u_{i_{k}}).
\end{align*}
For $\theta(t)=t\bs+(1-t)\bp$, with $\dot{\theta}(t)=\frac{d}{dt}\theta(t)=\bs-\bp$,
\[
\frac{1}{N}|H_{N}(\bs)-H_{N}(\bp)|\leq\frac{1}{N}\int_{0}^{1}|\dot{\theta}(t)\cdot\nabla H_{N}(\theta(t))|dt.
\]

Since $R_{s}(\dot{\theta}(t),\dot{\theta}(t))=R_{s}(\bs-\bp,\bs-\bp)$,
for $a:=\max_{s\in\S}\sqrt{R_{s}(\bs-\bp,\bs-\bp)}$ and some $u\in\overline{M}_{N}$,
$\dot{\theta}(t)=au$. Thus, the event from (\ref{eq:LipProb}) is
contained in the event that
\[
\frac{1}{N}\sup_{\bs,u\in\overline{M}_{N}}G(u,\bs)\leq L.
\]

By the Borell-TIS inequality, to prove the lemma it will be enough
to show that 
\begin{equation}
\frac{1}{N}\E\sup_{\bs,u\in\overline{M}_{N}}G(u,\bs)\leq C\label{eq:maxbd1}
\end{equation}
and that 
\begin{equation}
\frac{1}{N}\sup_{\bs,u\in\overline{M}_{N}}\E G(u,\bs)^{2}\leq C,\label{eq:maxbd2}
\end{equation}
for some $C>0$ independent of $N$.

Suppose $\bs$, $\bp$, $u$ and $v$ are points in $\overline{M}_{N}$.
From the inequality $(x_{1}+\cdots x_{k})^{2}\leq k(x_{1}^{2}+\cdots x_{k}^{2})$
and symmetry,
\begin{align*}
 & \E(G(u,\bs)-G(v,\bp))^{2}\\
 & =N\sum_{k=1}^{\infty}\sum_{i_{1},\dots,i_{k}=1}^{N}\Delta_{i_{1},\dots,i_{k}}^{2}[(u_{i_{1}}\cdots\sigma_{i_{k}}+\cdots\sigma_{i_{1}}\cdots u_{i_{k}})-(v_{i_{1}}\cdots\pi_{i_{k}}+\cdots\pi_{i_{1}}\cdots v_{i_{k}})]^{2}\\
 & \leq N\sum_{k=1}^{\infty}\sum_{i_{1},\dots,i_{k}=1}^{N}\Delta_{i_{1},\dots,i_{k}}^{2}k[(u_{i_{1}}\cdots\sigma_{i_{k}}-v_{i_{1}}\cdots\pi_{i_{k}})^{2}+(\sigma_{i_{1}}\cdots u_{i_{k}}-\cdots\pi_{i_{1}}\cdots v_{i_{k}})^{2}]\\
 & =N\sum_{k=1}^{\infty}\sum_{i_{1},\dots,i_{k}=1}^{N}\Delta_{i_{1},\dots,i_{k}}^{2}k^{2}(u_{i_{1}}\sigma_{i_{2}}\cdots\sigma_{i_{k}}-v_{i_{1}}\pi_{i_{2}}\cdots\pi_{i_{k}})^{2}\\
 & \leq N\sum_{k=1}^{\infty}\sum_{i_{1},\dots,i_{k}=1}^{N}\Delta_{i_{1},\dots,i_{k}}^{2}k^{3}\Big[(u_{i_{1}}-v_{i_{1}})^{2}\sigma_{i_{2}}^{2}\cdots\sigma_{i_{k}}^{2}\\
 & +\sum_{j=2}^{k}v_{i_{1}}^{2}\sigma_{i_{2}}^{2}\cdots\sigma_{i_{j-1}}^{2}(\sigma_{i_{j}}-\pi_{i_{j}})^{2}\pi_{i_{j+1}}^{2}\cdots\pi_{i_{k}}^{2}\Big],
\end{align*}
where for the last inequality we the fact that 
\begin{align*}
 & u_{i_{1}}\cdots\sigma_{i_{k}}-v_{i_{1}}\cdots\pi_{i_{k}}=(u_{i_{1}}-v_{i_{1}})\sigma_{i_{2}}\cdots\sigma_{i_{k}}\\
 & +v_{i_{1}}(\sigma_{i_{2}}-\pi_{i_{2}})\sigma_{i_{3}}\cdots\sigma_{i_{k}}+\cdots+v_{i_{1}}\pi_{i_{2}}\cdots\pi_{i_{k-1}}(\sigma_{i_{k}}-\pi_{i_{k}}).
\end{align*}
From this, one can check that
\begin{equation}
\E(G(u,\bs)-G(v,\bp))^{2}\leq N\Delta_{0}^{2}\sum_{s\in\S}\Big(R_{s}(u-v,u-v)+R_{s}(\bs-\bp,\bs-\bp)\Big),\label{eq:1spinc}
\end{equation}
for $\Delta_{0}^{2}:=\sum_{p\in P}\Delta_{p}^{2}|p|^{4}<\infty$.

Note that if we define 
\[
\tilde{G}(u,\bs)=\sum_{s\in\S}\sqrt{\frac{N}{N_{s}}}\sum_{i\in I_{s}}\Delta_{0}J_{i}\sigma_{i}+\sum_{s\in\S}\sqrt{\frac{N}{N_{s}}}\sum_{i\in I_{s}}\Delta_{0}\tilde{J}_{i}u_{i},
\]
where $J_{i}$ and $\tilde{J}_{i}$ are i.i.d. standard normal variables,
then $\E(\tilde{G}(u,\bs)-\tilde{G}(v,\bp))^{2}$ is equal to the
right-hand side of (\ref{eq:1spinc}). Hence, from the Sudakov\textendash Fernique
inequality, (see e.g. \cite[Theorem 2.2.3]{RFG})
\begin{equation}
\frac{1}{N}\E\sup_{\bs,u\in M_{N}}G(u,\bs)\leq\frac{1}{N}\E\sup_{\bs,u\in M_{N}}\tilde{G}(u,\bs)\leq2\Delta_{0}\sum_{s\in\S}\sqrt{N_{s}/N},\label{eq:EsupG}
\end{equation}
which proves (\ref{eq:maxbd1}).

A similar computation to the above gives, for any $\bs,\,u\in\overline{M}_{N}$,
\[
\frac{1}{N}\E G(u,\bs)^{2}\leq\sum_{k=1}^{\infty}\sum_{i_{1},\dots,i_{k}=1}^{N}\Delta_{i_{1},\dots,i_{k}}^{2}k^{2}u_{i_{1}}^{2}\sigma_{i_{2}}^{2}\cdots\sigma_{i_{k}}^{2}\leq\sum_{p\in P}\Delta_{p}^{2}|p|^{2},
\]
from which (\ref{eq:maxbd2}) follows.
\end{proof}
\bibliographystyle{plain}
\bibliography{master}

\end{document}